\documentclass[12pt]{amsart} 

\usepackage{amsmath,amsfonts,amssymb,amscd,graphics}

\def\dim{\mathop{\mathrm{dim}\;\!}\nolimits}

\def\zit{\mathbb{Z}}   
\def\nit{\mathbb{N}} 
\def\cit{\mathbb{C}} 
\def\B{\mathbb{B}}
\def\W{\mathbb{W}}

\def\D{\mathcal{D}}

\def\tr{\mathrm{tr}}
\def\gg{\mathbf{g}}
\def\gh{\mathbf{h}}
\def\gnm{\mathbf{n_-}}
\def\gnp{\mathbf{n_+}}

\def\gn{\mathbf{n}}

\def\bm{\nu_0}
\def\la{\langle}
\def\ra{\rangle}

\def\ts{d}

\def\ftext#1{{\let\thefootnote\relax\footnotetext{\noindent #1}}}
\def\bin#1#2{\left(\begin{smallmatrix}#1\\#2\end{smallmatrix}\right)}
\def\sw{\succ^{p}}

\newcommand{\al}{\alpha}

\renewcommand{\qed}{\hfill~~\mbox{$\square$}}

\newtheorem{theorem}{Theorem}[section]
\newtheorem{proposition}[theorem]{Proposition}

\newtheorem{lemma}[theorem]{Lemma}
\newtheorem{corollary}[theorem]{Corollary}

\begin{document}
\title[Hypergeometric solutions satisfy dynamical equations]{ Hypergeometric 
Solutions of Trigonometric KZ Equations satisfy
Dynamical Difference Equations}
\author{ Y. Markov$^\star$ and A. Varchenko$^{\star,1}$ }
\date{}
\maketitle
\begin{center}
{\it $^{\star}$  Department of Mathematics, University of North Carolina, Chapel Hill,\\
NC 27599 -- 3250, USA\\
{\{yavmar, anv\} @ email.unc.edu}}
\end{center}

\begin{abstract}
The trigonometric KZ equations associated to a Lie algebra $\gg$ depend on a parameter $\lambda\in\gh$
where $\gh\subset\gg$ is a Cartan subalgebra. A system of dynamical difference equations with respect to $\lambda$
compatible with the KZ equations is introduced in \cite{TV}. We prove that the standard hypergeometric 
solutions of the trigonometric KZ equations associated to $sl_N$ also satisfy the dynamical difference equations.
\end{abstract}\hfill\\

{\bf Mathematics subject classification:}{ Primary: 35Q40; Secondary 17B10}\hfill\\


\ftext{ \hspace{-0.6cm}
$^1$ Supported in part by NSF grant DMS-9801582.}

\section{Introduction}
Consider the simplest of hypergeometric integrals,
\begin{equation*}
 I(z,a,b)=\int_0^z t^{a-1}(t-z)^{b-1}\, dt.
\end{equation*}
It satisfies the differential equation
\begin{equation}\label{simpleKZ}
\frac{dI}{dz}(z,a,b)=\frac{a+b-1}{z}I(z,a,b),
\end{equation}
the difference equation
\begin{equation}\label{simpleDyn}
I(z,a+1,b)=\frac{az}{a+b}I(z,a,b),
\end{equation}
and a similar difference equation with respect to $b$. The equations (\ref{simpleKZ}) and 
(\ref{simpleDyn}) are clearly compatible. More general hypergeometric integrals appear in conformal field 
theory as integral representations for conformal blocks, see \cite{Ch,CF,DF,Mat,SV1,V1}. 
It is known that the integrals 
satisfy the KZ differential equations. The KZ equations are generalizations of (\ref{simpleKZ}). In \cite{TV}
a system of difference equations, generalizing equation (\ref{simpleDyn}), is proposed. The system, 
proposed in \cite{TV}, is compatible with the KZ equations. The difference equations were called the dynamical
equations. Both the KZ differential equations and the dynamical difference equations are associated to a 
given Lie algebra. It was conjectured in \cite{TV} that the hypergeometric integrals, which satisfy the KZ equations,
also satisfy the dynamical equations.

In this paper we prove that the hypergeometric integrals solving the (trigonometric) KZ equations 
associated to $sl_N$ also satisfy the the dynamical equations.

  The trigonometric KZ equations have a rational limit called the (standard) rational KZ equations. 
Under this limiting procedure the dynamical difference equations turn into a system of (dynamical) differential 
equations compatible with the rational KZ equations. The dynamical differential equations were introduced 
and studied in \cite{FMTV}, see also \cite{T-L}. In \cite{FMTV} it was shown that the  hypergeometric 
solutions of the rational KZ equations also satisfy the dynamical differential equations.

The paper is organized as follows. In Sections~\ref{TKZsec} and \ref{sec2} we introduce notation and
define the main objects of our study: the trigonometric and  rational KZ equations, 
the dynamical equations for the Lie algebra $sl_N$. 

In Section~\ref{HGsol}, following \cite{Mat}, we present a construction of 
hypergeometric solutions of the KZ equations with values in a tensor product of highest weight 
$sl_N$-modules.  

  The main result of the paper is Theorem~\ref{main_sln}. We give new formulae for hypergeometric solutions 
related to special (normal) orders on the set of positive roots of $sl_N$ in Sections~\ref{spec_ord} and 
\ref{Pf.5.2}, and new formulae for the dynamical 
difference equations in terms of the Shapovalov form in Section~\ref{add_sect}. Both
results are used in the proof of  Theorem~\ref{main_sln}, given in Section~\ref{mainproof}. 

  In Appendices A and B we adapt a theorem from \cite{Mat} and a theorem from \cite{EFK}, respectively,
to our setting. In Appendix C we give explicit formulae illustrating the main objects of our study in the 
case of the Lie algebra $sl_3$.

We thank P. Etingof, K. Styrkas, and V. Tarasov for many valuable discussions.

\section{Rational and trigonometric KZ equations}\label{TKZsec}
\subsection{Preliminaries}\label{prelim}
 Let $\gg$ be a simple complex Lie algebra with  root space decomposition $\gg =
\gh\oplus(\oplus_{\alpha\in\Sigma}\cit e_{\alpha})$
 where $\Sigma\subset\gh^*$ is the set of  roots. Fix a system of
simple roots $\alpha_1,\ldots ,\alpha_r$. Let $\Sigma_{\pm}$ be the set of positive
(negative) roots. Let $\gn_{\pm}=\oplus_{\alpha\in \Sigma_{\pm}}\gg_\alpha$. Then
$\gg=\gn_+\oplus\gh\oplus\gn_-$.

Let $(\,,\,)$ be an invariant bilinear form on $\gg$. The form gives rise to a 
natural identification $\gh\to\gh^*$. We use this identification
and make no distinction between $\gh$ and $\gh^*$.
This identification allows us to define a scalar product on $\gh^*$.
We  use the same notation $(\,,\,)$ for the pairing $\gh\otimes \gh^*\to\cit$.

We use the  notation: $Q=\oplus_{k=1}^r\zit\alpha_k$ - root lattice;
$Q^+=\oplus_{k=1}^r\zit_{\ge 0}\alpha_k$;
$Q^\vee=\oplus_{k=1}^r\zit\alpha_k^\vee$ - dual root lattice,
where $\alpha^\vee=2\alpha/(\alpha,\alpha)$;
$P=\{\lambda\in\gh\,|\, (\lambda,\alpha^\vee_k)\in\zit\}$ - weight lattice;
$P^+=\{\lambda\in\gh\,|\, (\lambda,\alpha^\vee_k)\in\zit_{\ge 0}\}$ - 
cone of dominant integral weights;
$\omega_k\in P^+$ - fundamental weights: $(\omega_k,\alpha^\vee_l)=\delta_{k,l}$;
$\rho={\frac12}\sum_{\alpha\in\Sigma_+}\alpha=\sum_{k=1}^r\omega_k$;
$P^\vee=\oplus_{k=1}^r\zit\omega^\vee_k$ - dual weight lattice, where $\omega^\vee_k$
-dual fundamental weights: $(\omega^\vee_k,\alpha_l)=\delta_{k,l}$.

Define a partial order on $\gh$  putting $\mu<\lambda$ if $\lambda-\mu\in Q^+$.

For $\alpha\in\Sigma$ choose generators $e_\alpha\in\gg_\alpha$ so that 
$(e_\alpha,e_{-\alpha})=1$. For any $\alpha$, the triple 
$$ H_\al=\al^\vee, \qquad E_\al=\frac{2}{(\al,\al)}e_\al,\qquad F_\al=e_{-\al}$$
forms an $sl_2$-subalgebra in $\gg$, $[H_\al,E_\al]=2E_\al,\, [H_\al,F_\al]=-2F_\al,\, 
[E_\al,F_\al]=H_\al$.

 The Chevalley involution $\tau$ of $\gg$ is defined by
$E_{\alpha_k}\mapsto -F_{\alpha_k}$, $F_{\alpha_k}\mapsto
 -E_{\alpha_k}$, $H_{\alpha_k}\mapsto -H_{\alpha_k}$, $ k=1,\ldots,r$.
 The antipode map $A$ of $\gg$ is defined by 
$g\mapsto -g$ for $g\in \{E_{\alpha_k},F_{\alpha_k}, H_{\alpha_k}\}_{k=1}^{r}$.

 Let $U(\gg)$ be the universal enveloping algebra of $\gg$.
The Chevalley involution, $\tau$, extends to an involutive automorphism of $U(\gg)$ which
permutes $U(\gnp)$ and $U(\gnm)$.  The antipode map, $A$,
extends to an involutive anti-automorphism of $U(\gg)$ which preserves $U(\gnp)$ and $U(\gnm)$.

Let $s_k:\gh\to\gh$ denote 
the simple reflection $s_k(\lambda)=\lambda-(\alpha_k^\vee,\lambda)\alpha_k$ for all $\lambda\in\gh$, 
and let $\W$ be the Weyl group, generated by $s_1,...,s_r$.
For an element $w\in \W$,  denote $l(w)$ 
the length of the minimal (reduced) presentation of $w$ as a product of generators
$s_1,...,s_r$.

 For any dual fundamental weight $\omega_k^\vee$ define
an element $\omega_{[k]}=\omega_0\omega_0^k\in \W$ where $\omega_0$ (respectively, $\omega_0^k$) 
is the longest element in $\W$ ( respectively, in $\W^k$ generated by all simple reflections 
$s_l$ preserving $\omega_k^\vee$).

Let $\lambda\in\gh$ be a weight. Let $\cit_\lambda$
be the one-dimensional $(\gh\oplus\gnp)$-module such that $\cit_\lambda=\cit v_\lambda$ with
$h\,v_\lambda=\lambda(h)v_\lambda$ for any $h\in\gh$ and
$\gnp\,v_\lambda=0$. The Verma module with highest weight
$\lambda$ is the induced module $M_\lambda=\mathrm{Ind}_{(\gh\oplus\gnp)}^\gg \cit_\lambda$.
$M_\lambda$ is a free $U(\gnm)$- module and can be identified with $U(\gnm)$ as a linear space
by the map $U(\gnm)\rightarrow M_\lambda$, $u\mapsto u\,v_\lambda$ for any $u\in U(\gnm)$.
A highest weight module of weight $\lambda$ is a quotient module of the
Verma module with highest weight $\lambda$.
All the highest weight modules in this paper are equipped with a distinguished generator, 
the highest vector.

 Let $V_\lambda$ be a highest weight $\gg$-module with highest weight
$\lambda$ and  highest weight vector $v_\lambda$. 
 We have a weight decomposition $V_\lambda=\otimes_{\nu\leq\lambda}V_\lambda[\nu]$.
Define $V^*_\lambda=\otimes_{\nu\leq\lambda}V_\lambda[\nu]^*$ the restricted dual module to $V_\lambda$ with the 
$\gg$-action $\la g\phi, a\ra=-\la \phi, g a\ra$ 
for every $g\in\gg$, $a\in V_\lambda$, $\phi\in V^*_\lambda$. Then $V^*_\lambda$ is a lowest weight 
$\gg$-module with the homogeneous lowest weight vector $v^*_\lambda$, such that $\la v_\lambda^*,v_\lambda\ra =1$.

 The Shapovalov form $S_\lambda(\,.\,,\,.\,)$ on $V_\lambda$  is 
the unique symmetric bilinear form such that
$$ S_\lambda(v_\lambda,v_\lambda)=1,\quad S_\lambda(g u,v)=S_\lambda(u,(A\circ\tau)(g) v),
\quad \forall\, u,v\in V_\lambda,\, \forall\, g\in\{E_{\alpha_k},F_{\alpha_k}\}_{k=1}^r.$$
The form $S_\lambda(\,.\,,\,.\,)$ is non-degenerate if and only if $V_\lambda$ is irreducible. 
In particular, it is non-degenerate for generic values of $\lambda$.


\subsection{ KZ equations.}
Let $\{x_k\}$ be any
orthonormal basis of the Cartan subalgebra $\gh$. Set
$$\Omega^0=\frac12\sum_k x_k\otimes x_k,\qquad
\Omega^+=\Omega^0+\sum_{\alpha\in\Sigma_+}e_{\alpha}\otimes e_{-\alpha},\qquad
\Omega^-=\Omega^0+\sum_{\alpha\in\Sigma_+}e_{-\alpha}\otimes e_{\alpha}.$$
Define the Casimir operator $\Omega$ and the trigonometric R-matrix $r(z)$ by
$$ \Omega=\Omega^+ +\Omega^-,\qquad 
r(z)=\frac{\Omega^+z +\Omega^-}{z-1}.$$

Let $V'=V_1 \otimes\cdots\otimes V_{n+1}$, where $V_j$ is a 
$\gg$--module. The {\it rational KZ operators},  $\nabla_{KZ,i}(\kappa)$,
acting on a function $u(z_1,\ldots,z_{n+1})$ of $n+1$ complex 
variables with values in $V'$ are
\begin{equation}\label{KZop}
 \nabla_{KZ, i}(\kappa)= \kappa\frac{\partial}{\partial z_i}- \sum_{j, j\ne i}
 \frac{\Omega^{(ij)}}{z_i - z_j}, \qquad\qquad i=1,\ldots,n+1,
\end{equation}
where $\kappa$ is a complex parameter. The {\it rational KZ equations} are
\begin{equation}\label{KZ}
 \nabla_{KZ, i}(\kappa)u(z_1,\ldots,z_{n+1})=0,\qquad i=1,\ldots,n+1. 
\end{equation}
The rational KZ equations are compatible, $[\nabla_{KZ, i},\nabla_{KZ, j}]=0$.
The KZ operators commute with the $\gg$ action on $V'$. Thus, the KZ operators 
preserve every subspace of $V'$ consisting of all singular vectors of a given weight. \\

  Let  $V=V_1 \otimes\cdots\otimes V_{n}$, where $V_j$ is a
 $\gg$--module. The {\it trigonometric KZ operators},
$\nabla_{i}(\kappa,\lambda)$, with parameters $\kappa \in \cit$ and $\lambda\in \gh$
 acting on a function $v(z_1,\ldots,z_{n},\lambda)$ of $n$ complex 
variables with values in $V$ are
\begin{equation}\label{TKZop}
\nabla_i(\kappa,\lambda)=\kappa z_i\frac{\partial}{\partial z_i}-
\lambda^{(i)} - \sum_{j, j\ne i}
 r(z_i/z_j)^{(ij)}, \qquad\qquad i=1,\ldots,n.
\end{equation}
 The {\it trigonometric KZ equations} are
\begin{equation}\label{TKZ}
 \nabla_{ i}(\kappa,\lambda)v(z_1,\ldots,z_{n},\lambda)=0,\qquad i=1,\ldots,n. 
\end{equation}
 The trigonometric KZ equations are compatible, $[\nabla_{ i},\nabla_{j}]=0$.
The trigonometric KZ operators commute with the $\gh$ action on $V$. Thus, the trigonometric 
KZ operators preserve every weight subspace of $V$.

\subsection{A relation between rational KZ equations and 
 trigonometric KZ equations} For every $j=1,\ldots,n+1$, let $V_j$ be
a highest weight $\gg$-module with highest weight $\Lambda_j$ and
highest weight vector $v_j$. Set $V'=V_1 \otimes\cdots\otimes V_{n+1}$
and  $V=V_1 \otimes\cdots\otimes V_{n}$. 

Let $V_{n+1}^*$ be the dual module to $V_{n+1}$ with the homogeneous lowest weight vector $v^*_{n+1}$.
 Define a multi-linear map
$|v_{n+1}^*\ra: V'\rightarrow V$ by  
$y_1\otimes\cdots\otimes y_n\otimes y_{n+1}|v_{n+1}^*\ra
=\la v_{n+1}^*, y_{n+1}\ra y_1\otimes\cdots\otimes y_n$, 
for any $y_1\otimes\cdots\otimes y_n\otimes y_{n+1}\in V'$.
  
  The following well known fact, see \cite{EFK} for example, describes the transition 
from rational KZ equations to trigonometric KZ equations.
\begin{proposition}\label{KZ2tKZ}
Fix a weight subspace $V'[\nu']\subset V'$, $\nu'=\sum_{j=1}^{n+1}\Lambda_j-\nu_0$, where  
$\bm\in Q_+$.  Let $u:\cit^{n+1}\rightarrow V'$ be a solution of the rational KZ 
equations with parameter $\kappa\in\cit$ taking values in the subspace of  $V'[\nu']$ 
consisting of all singular vectors. 
Set $\nu=\sum_{j=1}^{n}\Lambda_j-\nu_0$.

Then $v(z_1,\cdots,z_n)=u(z_1,\cdots,z_n,0)|v_{n+1}^*\ra
\prod_{i=1}^n z_i^\frac{(\Lambda_i,\Lambda_i+2\rho)}{2\kappa}$  is a solution of the 
trigonometric KZ equations with values in the weight subspace $V[\nu]\subset V$ with parameter
$\lambda = \Lambda_{n+1} +\rho +\frac12\nu\in\gh$ and the same parameter $\kappa\in\cit$.
\end{proposition}
A proof is given in Appendix B. \qed
  
\section{Dynamical difference equations for $\gg=sl_N$, \cite{TV}.}\label{sec2}
 
\subsection{Operators $\B^\alpha_V(\lambda)$, $\B_{\omega,V}(\lambda)$.}\label{def_B}
 Let $\gg$ be a 
simple complex Lie algebra.
Fix a root $\alpha\in\Sigma_+$. Consider the $sl_2$ subalgebra of $\gg$ with generators
$H=H_\alpha$, $E=E_\alpha$, $F=F_\alpha$. For $t\in\cit$, introduce
\begin{equation}\label{p-oper}
 p(t;H,E,F)=\sum_{k=0}^\infty F^kE^k\frac1{k!}\prod_{j=0}^{k-1}\frac1{(t-H-j)}.
\end{equation}
The series $p(t,H,E,F)$ is an element of a suitable completion of $U(sl_2)$.

  Let $V=V_1 \otimes\cdots\otimes V_{n}$ be a tensor product of highest weight 
$\gg$--modules.  For $\alpha\in\Sigma$, $\lambda\in\gh$, the linear operator 
$\B^\alpha_V(\lambda):V\rightarrow V$ is defined by the rule:
for any $\nu\in\gh$, and any $v\in V[\nu]$, we have
$$\B^\alpha_V(\lambda)v=p((\lambda+\frac12\nu,\alpha^\vee)-1;H_{\alpha},E_\alpha,F_\alpha)v.$$

  For every simple reflection $s_k\in \W$ define $\B_{s_k,V}(\lambda):V\rightarrow V$ by 
$\B_{s_k,V}(\lambda)=\B^{\alpha_k}_V(\lambda)$.
 For every $w\in\W$, such that $l(w)>1$, the operator 
$ \B_{w,V}(\lambda)$ is defined by the rule.
If $w=w_1w_2$, where $w_1,w_2\in\W$ and $l(w_1w_2)=l(w_1)+l(w_2)$, then 
$$ \B_{w,V}(\lambda)=w_2^{-1}(\B_{w_1,V}(w_2(\lambda)))\B_{w_2,V}(\lambda).$$ 

\subsection{The Lie algebra $sl_N$.}  Let $\{e_{k,l}\}_{k,l}$,
 $k,l=1,\ldots,N$, be the standard generators of the Lie algebra $gl_N$,
$[ e_{k,l}\,, \, e_{k',l'}]\,=\,\delta_{l,k'}\,e_{k,l'}\,-\,\delta_{k,l'}\,e_{l,k'}\,.$
The Lie algebra $sl_N$ is the Lie subalgebra of $gl_N$ such that  
$sl_N=\gn_+\oplus\gh\oplus\gn_-$ where
$$
\gn_+=\oplus_{1\leq k < l\leq N}\cit\,e_{k,l}\,,\qquad
\gn_-=\oplus_{1\leq k < l\leq N}\cit\,e_{l,k}\,, 
$$
and $\gh=\{ \lambda=\sum_{k=1}^N\lambda_ke_{k,k}\,|\,\lambda_k\in\cit,\,\,
\sum_{k=1}^N\lambda_k=0\}$.

The invariant scalar product on $sl_N$ is defined by 
$(e_{k,l}, e_{k',l'})=\delta_{k,l'}\delta_{l,k'}$.

The roots of $sl_N$ are $\alpha_{k,l}=e_{k,k}-e_{l,l}$, $k\neq l$. The positive roots are 
$\{\alpha_{k,l}\}_{k<l}$. The simple roots are $\alpha_k=e_{k,k}-e_{k+1,k+1}$, 
$k=1,...,N-1$.
For a positive root $\alpha_{k,l}$, the elements $H_{\alpha_{k,l}}=e_{k,k}-e_{l,l}$,
$E_{\alpha_{k,l}}=e_{k,l}$, $F_{\alpha_{k,l}}=e_{l,k}$ generate 
the $sl_2$ subalgebra associated with $\alpha_{k,l}$.

We have $\alpha^\vee=\alpha$ for any root.

 For any $k<l$, we have
\begin{equation}\label{e_prod_f}
F_{\alpha_{k,l}}=[F_{\alpha_{l-1}},[\ldots,[F_{\alpha_{k+1}},F_{\alpha_k}]\ldots]].
\end{equation}

Define the standard linear order 
on the set of positive roots,  $\alpha_{k,l}\succ \alpha_{k',l'}$ if and only if $l>l'$, 
or $l=l'$ and $k>k'$.  

The Weyl group
$\W$ is the symmetric group $S^N$ permuting  coordinates of $\lambda\in\gh$.
The (dual) fundamental weights are 
$\omega_k=\omega_k^\vee=\sum_{h=1}^k(1-\frac{k}{N})e_{h,h}
-\sum_{h=k+1}^N\frac{k}{N}e_{h,h}$,
 and  the permutations $w_{[k]}^{-1}\in S^N$ have the form
$w_{[k]}^{-1}=\left( {}^1_{k+1}\,{}^{2}_{k+2}\,{}^{...}_{...}\,{}^{N-k}_{N}\,{}^{N-k+1}_{1}\,
{}^{...}_{...}\,{}^{N}_{k} \right)$, $k=1,\ldots,N-1$.

\subsection{Dynamical difference equations}
The {\it dynamical difference equations} on a $V$-valued function $v(z_1,\ldots,z_n,\lambda)$
for $sl_N$ are
\begin{equation}\label{DynEq}
v(z_1,\ldots,z_n,\lambda+\kappa\omega_k^\vee)=
K_k(z_1,\ldots,z_n,\lambda)v(z_1,\ldots,z_n,\lambda), \qquad k=1,\ldots, N-1,
\end{equation}
where
$$K_k(z_1,\ldots,z_n,\lambda)=\prod_{j=1}^nz_j^{(\omega_k^\vee)^{(j)}}
\B_{\omega_{[k]},V}(\lambda).$$ 

The operators $K_k(z,\lambda)$ preserve the weight decomposition of $V$.
\begin{theorem}[Theorem 17 in \cite{TV}] \it The dynamical equations (\ref{DynEq}) together 
with the trigonometric KZ equations (\ref{TKZ}) form a compatible system of
equations. Namely, 
\begin{align}
&[\nabla_i(\kappa,\lambda),\nabla_{j}(\kappa,\lambda)]=0, \qquad
\nabla_j(\kappa,\lambda+\kappa\omega_k^\vee)K_k(z,\lambda)=
K_k(z,\lambda)\nabla_j(\kappa,\lambda),\notag\\
&K_k(z,\lambda+\kappa\omega_{l}^\vee)K_{l}(z,\lambda)=
K_{l}(z,\lambda+\kappa\omega_{k}^\vee)K_{k}(z,\lambda)  
\end{align}
 for all $i,j=1,\ldots,n$, and $k,l=1,\ldots,N-1$.
\end{theorem}


\section{Hypergeometric solutions of the trigonometric KZ equations for $sl_N$.}\label{HGsol}
We use the construction in \cite{Mat} of hypergeometric solutions of  
the rational $sl_N$ KZ equations in a tensor product of lowest weight modules.
 We modify this procedure in Appendix A to a construction of solutions of the rational KZ equations
 in a tensor product of 
highest weight modules, cf. \cite{SV1}. We present the result below. Then we  
use Proposition~\ref{KZ2tKZ} to give hypergeometric solutions of the trigonometric KZ 
equations.

There are different constructions of hypergeometric solutions of KZ equations, cf. \cite{Ch,CF,DF,Mat,SV1,V1}.
One should expect, that all of them give the same result, but this was never checked off as far as we know.

 Let $V'=V_1\otimes\cdots\otimes V_{n+1}$, where  $V_j$ is a highest 
weight $sl_N$ module with  highest weight $\Lambda_j$ and highest weight vector $v_j$.
Fix a weight subspace  $V'[\nu']\subset V'$, $\nu'=\sum_{j=1}^{n+1}\Lambda_j-\sum_{k=1}^{N-1}m_k\alpha_k$,
where $\sum_{k=1}^{N-1} m_k\alpha_k\in Q^+$.
Set $m=\sum_{k=1}^{N-1} m_k$ and $\bm=\sum_{k=1}^{N-1}m_k\alpha_k$.

{\bf The function $\Phi'$.} Consider complex spaces $\cit^{n+1}$ with coordinates $z_1,\ldots,z_{n+1}$,
and $\cit^m$ with coordinates $t_k^{(\ts)}$, $k=1,\ldots, N-1$, $\ts\in S_k=\{1,\ldots,m_k\}$. 
Fix an order $\ll$  on the set of coordinates in $\cit^m$, 
$(k,\ts) \ll (k',\ts')$ if and only if $k<k'$, or $k=k'$ and $\ts<\ts'$. Define a multi-valued function
$\Phi':\cit^{n+1}_{z'}\times\cit^m_t\rightarrow\cit$
\begin{equation}
\Phi'(z',t)=\prod_{i<j}(z_i-z_j)^{(\Lambda_i,\Lambda_j)}
\prod_{(k,\ts),j}(t_k^{(\ts)}-z_j)^{-(\alpha_k,\Lambda_j)}
\prod_{(k,\ts)\ll (l,\ts')}(t_k^{(\ts)}-t_l^{(\ts')})^{(\alpha_{k},\alpha_{l})}.
\end{equation}

{\bf The standard PBW-basis.}
 Any order on the set $\Sigma_+$ of positive roots of $sl_N$  induces an order on the set
$\{-e_{l,k}\}_{k<l}$ which is a basis of $\gnm$, $-e_{l,k}$ succeeds $ -e_{l',k'}$ if and only if
$\alpha_{k,l}$ succeeds $\alpha_{k',l'}$. The standard order $\succ$ on $\Sigma_+$ induces the standard 
order $\succ$ on $\{-e_{l,k}\}_{k<l}$,
 $-e_{l,k}\succ -e_{l',k'}$ if and only if $l>l'$, or $l=l'$ and $k>k'$.
The corresponding (standard) PBW-basis of $U(\gnm)$ is 
$$\left\{ F_{I_0} = (-1)^{\sum i_{l,k}}\frac{e_{N,N-1}^{i_{N,N-1}}}{i_{N,N-1}!}\cdots
\frac{e_{2,1}^{i_{2,1}}}{i_{2,1}!}\right\},$$
where $I_0=\{i_{l,k}\}_{l>k}$  runs over all sequences of non-negative integers. 

Let $F_{I_1},\ldots,F_{I_{n+1}}$ be elements of the standard PBW-basis,   
$I_j=\{i^j_{2,1},\ldots,i^j_{N,N-1}\}$. Set $I=(I_1,\ldots,I_{n+1})$.
The corresponding ``monomial'' vector  $F_Iv=F_{I_1}v_1\otimes\cdots\otimes F_{I_{n+1}}v_{n+1}$ 
lies in $V'[\nu']$ if 
\begin{equation}\label{pmn} 
\sum_{j=1}^{n+1}\sum_{k=1}^h\sum_{l=h+1}^N i^j_{l,k}=m_h,\quad\mbox{ for all } h=1,\ldots,N-1.
\end{equation}
Denote $P(\bm,n+1)$ the set of all indices $I$ corresponding to monomial vectors in $V'[\nu']$. 
The set $\{ F_Iv\}_{I\in P(\bm,n+1)}$ forms a basis of $V'[\nu']$ provided 
the tensor factors of $V'$ are Verma modules, and generates $V'[\nu']$ otherwise.

{\bf The rational function $\phi'$.} 
For  $I\in P(\bm,n+1)$ and  $h=1,\ldots,N-1$, define two index sets.
\begin{align}
S(I)&=\{ (j,k,l,q)\,|\, 1\leq j\leq n+1,\quad 1\leq k<l\leq N,\quad 1\leq q\leq i^j_{l,k}\}\notag\\
S_h(I)&= \{ s=(j,k,l,q)\in S(I)\,|\, k\leq h<l\}.  
\end{align}

Condition (\ref{pmn}) implies $|S_h(I)|=m_h$ for all $h=1,\ldots,N-1$. 

For every $h$ fix a bijection $\beta_h(I):S_h(I)\rightarrow S_h$. 
  For $s=(j,k,l,q)\in S(I)$, define rational functions 
\begin{equation}\label{eq17.5}
f^{(s)}=\prod_{h=k}^{l-2}\frac1{t_h^{(\beta_h(I)(s))}-t_{h+1}^{(\beta_{h+1}(I)(s))}},\qquad
\phi^{(s)}=f^{(s)}\frac1{t_{l-1}^{(\beta_{l-1}(I)(s))}-z_j}.
\end{equation}
Set
\begin{equation}
  \phi(I)=\prod_{s\in S(I)}\phi^{(s)}, \qquad \phi'(z',t)=\sum_{I\in P(\bm,n+1)}\phi(I)F_Iv.
\end{equation}

{\bf Examples.} Let $\gg=sl_3$. Then $\{-e_{3,2},-e_{3,1},-e_{2,1}\}$ is a basis of $\gnm$.
Let $n=0$.
\begin{align}
\notag &\mbox{ (a) } \nu'=\Lambda_1-\alpha_1-\alpha_2,\quad
 \phi'=\frac1{(t_1^{(1)}-z_1)(t_2^{(1)}-z_1)}e_{3,2}e_{2,1}v_1-
\frac1{(t_1^{(1)}-t_2^{(1)})(t_2^{(1)}-z_1)}e_{3,1}v_1.\\
\notag & \mbox{ (b) } \nu'=\Lambda_1-2\alpha_1-\alpha_2,\\
\notag & \phi'=
 \frac1{(t_1^{(1)}-t_2^{(1)})(t_2^{(1)}-z_1)(t_1^{(2)}-z_1)}e_{3,1}e_{2,1}v_1-
\frac{1}{(t_1^{(1)}-z_1)(t_1^{(2)}-z_1)(t_2^{(1)}-z_1)}e_{3,2}\frac{e_{2,1}^2}{2}v_1.
\end{align}

{\bf The integrals.} Consider the integral with values in $V'[\nu']$
$$u(z')=\int_{\gamma(z')}\Phi'(z',t)^{\frac1{\kappa}}\phi'(z',t)\,dt,$$ 
 where $dt=\prod_{(k,\ts)}dt^{(\ts)}_k$ and $\gamma(z')$ in $\{ z'\}\times\cit^m_t$ is
a horizontal family of $m$-dimensional cycles of the twisted homology defined 
by the multi-valued function $(\Phi')^{\frac1{\kappa}}$, see \cite{Mat}, \cite{SV1}, \cite{V1}.

For a positive integer $m_0$,  let $\Sigma_{m_0}$ be the symmetric group on $m_0$ elements.
The group $\Sigma(m_1,\ldots,m_{N-1})=\Sigma_{m_1}\times\ldots\times\Sigma_{m_{N-1}}$ acts on points
$t=\{t_k^{(\ts)}\,|\,k=1,\ldots,N-1,\,\ts\in S_k\}$ of $\cit^m$ permuting coordinates in each group
$\{t_k^{(\ts)}\,|\,\ts\in S_k\}$.
 Denote $\D(z')$ the union of hyperplanes $\displaystyle{\cup_{(k,\ts),j}}\{t_k^{(\ts)}=z_j\}\cup
\displaystyle{\cup_{(k,\ts),(k+1,\ts')}}\{t_k^{(\ts)}=t_{k+1}^{(\ts')}\}$ in $\cit^m_t$. 
 We always make the following assumption on $\gamma(z')$.

{\bf Assumption, \cite{Mat}.} {\it For any rational function $\phi$ with poles in $\D(z')$ and any
permutation\\ $\sigma\in\Sigma(m_1,\ldots,m_{N-1})$ we have 
$\int_{\gamma(z')}\Phi'(z',t)^{\frac1{\kappa}}\phi'(z',t)\,dt= 
\int_{\gamma(z')}\Phi'(z',\sigma t)^{\frac1{\kappa}}\phi'(z',\sigma t)\,dt$.}

\begin{theorem}[Theorem~2.4 in \cite{Mat}, and Corollary~\ref{cor_mat} of the present paper]\label{t_mat}
The function
  $$u(z')=\int_{\gamma(z')}\Phi'(z',t)^{\frac1{\kappa}}\phi'(z',t)\,dt$$
takes values in the subspace of  singular vectors in $V'[\nu']$ and satisfies the
rational KZ equations with parameter $\kappa\in\cit$.
\end{theorem}
The function $u$ is called a hypergeometric solution of the rational KZ equation. Different solutions 
correspond to different choices of the horizontal family $\gamma(z')$.

{\bf Remark.} The assumption on the cycles of integration implies that the function $u(z')$ 
does not depend on the choice of bijections $\{\beta_h(I)\}$.\\

 Let $V=V_1\otimes\cdots\otimes V_{n}$.
Fix a weight subspace  $V[\nu]\subset V$, $\nu=\sum_{j=1}^{n}\Lambda_j-\nu_0$.
 Let $z=(z_1,\ldots,z_n)\in \cit^n$. Define a function
\begin{align}
\notag\Phi(z,t;\lambda)&=\prod_{i<j}(z_i-z_j)^{(\Lambda_i,\Lambda_j)}
\prod_{(k,\ts),j}(t_k^{(\ts)}-z_j)^{-(\alpha_{k},\Lambda_j)}
\prod_{(k,\ts)<(l,\ts')}(t_k^{(\ts)}-t_l^{(\ts')})^{(\alpha_{k},\alpha_{l})}\times\\
&\times\prod_{(k,\ts)} (t_k^{(\ts)})^{-(\alpha_{k},\lambda-\rho-\nu/2)}\times
\prod_{i=1}^nz_i^{(\Lambda_i,\lambda-\nu/2+\Lambda_i/2)}.
\end{align}
Define a rational function $\phi(z,t)$ with values in the weight subspace $V[\nu]$,
$$\phi(z,t)=\sum_{I\in P(\bm,n)}\phi(I)F_Iv.$$
Consider the integral with values in $V[\nu]$, 
$\int_{\gamma(z)}\Phi(z,t;\lambda)^{1/k}\phi(z,t)\,dt$, where
$\gamma(z)$ in $\{ z\}\times\cit^m$ is
a horizontal family of $m$-dimensional cycles of the twisted homology defined by the 
multi-valued function $(\Phi)^{\frac1{\kappa}}$. 
Consider the union of hyperplanes $\D=\D(z_1,\ldots,z_n,0)$ in $\cit^m_t$. 
We always make the following assumption on $\gamma(z)$.  

{\bf Assumption.} {\it For any rational function $\phi$ with poles in $\D$ 
and any permutation $\sigma\in\Sigma(m_1,\ldots,m_{N-1})$ we have 
$\int_{\gamma(z)}\Phi(z,t;\lambda)^{1/\kappa}\phi(z,t)\,dt= 
\int_{\gamma(z)}\Phi(z,\sigma t;\lambda)^{1/\kappa}\phi(z,\sigma t)\,dt$.}

\begin{corollary}\label{mat_trig} The function
$$v(z;\lambda)=\int_{\gamma(z)}\Phi(z,t;\lambda)^{1/k}\phi(z,t)\,dt$$ 
takes values in the weight space $V[\nu]$ and satisfies
the trigonometric KZ equations  with parameters $\kappa\in\cit$ and
$\lambda\in\gh$. 
\end{corollary}
{\it Proof.} The statement of the corollary follows from Theorem~\ref{t_mat} and 
Proposition~\ref{KZ2tKZ}.\qed

The function $v$ is called a hypergeometric solution of the trigonometric KZ equation.

\section{ The main result}
 For any $k\in\{1,\ldots,N-1\}$ we have 
$$\Phi^{\frac{1}{\kappa}}(z,t;\lambda+\kappa\omega_k^\vee)=\left(
\prod_{j=1}^nz_j^{(\Lambda_j,\omega_k^\vee)}\prod_{d=1}^{m_k}\frac1{t_k^{(d)}}\right)
\Phi^{\frac{1}{\kappa}}(z,t;\lambda).$$
The product $\prod_{d=1}^{m_k}(t_k^{(d)})^{-1}$ is a univalued function, which never vanishes. Therefore,
we can and will consider any  horizontal family of $m$-dimensional cycles, $\gamma(z)$ in $\{ z\}\times\cit^m$, 
of the twisted homology defined by the multi-valued function $(\Phi(z,t;\lambda))^{\frac1{\kappa}}$ as
a horizontal family of $m$-dimensional cycles of the twisted homology defined by the 
multi-valued function $(\Phi(z,t;\lambda+\kappa\omega_k^\vee))^{\frac1{\kappa}}$. This identification implies, that 
if we choose a hypergeometric solution of the trigonometric KZ equation, 
$v(z,\lambda)=\int_{\gamma(z)}\Phi(z,t;\lambda)^{1/k}\phi(z,t)\,dt$ for fixed $\lambda$,
we obtain a hypergeometric solutions 
$v(z,\lambda+\kappa\omega^\vee)=\int_{\gamma(z)}\Phi(z,t;\lambda+\kappa\omega^\vee)^{1/k}\phi(z,t)\,dt$ 
for any $\omega^\vee\in P^\vee$.

The following theorem is the main result of this paper.
\begin{theorem} \label{main_sln}
Let $V=V_1\otimes\cdots\otimes V_n$ be a tensor product of highest weight $sl_N$ modules. 
Let $v$ be a hypergeometric solution of the trigonometric KZ equations (\ref{TKZ}) 
with values in a weight subspace 
$V[\nu]\subset V,$ $\nu\in \gh$. Then, the function $v$ also satisfies the dynamical  
equations (\ref{DynEq}),
$$ v(z,\lambda+\kappa\omega_k^\vee)=K_k(z,\lambda)v(z,\lambda),\qquad
 k=1,\ldots,N-1. $$ 
\end{theorem}
  The theorem is proved in Section~\ref{mainproof}.

{\bf Example.} Let $\gg=sl_2$ and $n=1$. Denote $\alpha$ the positive root of $sl_2$, and let 
$V_1=L_p$, where $p$ is a positive integer and $L_p$
is the $(p+1)$-dimensional irreducible $sl_2$-module with highest weight $p\frac{\alpha}2$ 
and highest weight vector $v_p$. Consider a weight subspace
$L_p[\nu]$ of $L_p$, where $\nu=(p-2m)\frac{\alpha}2$. Let $v(z_1,\lambda)$ be a hypergeometric solution
of the trigonometric KZ equation with parameters $\lambda\in\gh$, $\kappa\in\cit$, which takes values in
$L_p[\nu]$. Up to a multiplicative constant, it has the form
\begin{equation}\label{ex_sl2}
 v(z,\lambda)= z_1^{\frac1{2\kappa}(p-2m)(\lambda,\alpha)}
I_m\left(-\frac1\kappa((\lambda,\alpha)-1-\frac{p-2m}2)+1, 
-\frac{p}{\kappa}, \frac1{\kappa}\right)e_{2,1}^mv_p,
\end{equation}
where $I_m(a,b,c)$ is the Selberg integral, see \cite{Me},
\begin{align}
I_m(a,b,c)&=\int_{0\leq t_1 < \cdots < t_m\leq 1}\left(\prod_{k=1}^n t_k^{a-1}(1-t_k)^{b-1}\right)
\prod_{1\leq k < l\leq m}(t_l-t_k)^{2c}\, dt_1\cdots dt_m\notag\\
&=\frac1{m!}\prod_{j=0}^{m-1}\frac{\Gamma(1+c+jc)\Gamma(a+jc)\Gamma(b+jc)}{\Gamma(1+c)\Gamma(a+b+(m+j-1)c)}.\notag
\end{align}
The dynamical equation (\ref{DynEq}) reduces to the following equation satisfied by the Selberg integral
$$I_m(a+1,b,c)=\left(\prod_{k=1}^m\frac{a+c(m-k)}{a+b+c(2m-k-1)}\right)I_m(a,b,c).$$

{\bf Application to determinants.} Let $V$ be a finite dimensional $sl_N$ - module, and
$V[\nu]$ a weight subspace. For a positive root $\alpha$ fix the $sl_2$ subalgebra in $sl_N$ generated by 
$H_\alpha, E_\alpha, F_\alpha$. Consider V as an $sl_2$ module. Let $V[\nu]_\alpha\subset V$ be the $sl_2$ - 
submodule generated by $V[\nu]$. Let 
$V[\nu]_\alpha=\oplus_{m\in\zit_{\geq}}W_m^\alpha\otimes L_{\nu+m\alpha}$ be the decomposition into 
irreducible $sl_2$ - modules, where $L_{\nu+m\alpha}$ is the irreducible module with highest weight
$\nu+m\alpha$ and $W_m^\alpha$ is the multiplicity space. Set $d_m^\alpha=\dim W_m^\alpha$ and 
$$X_{\alpha,V[\nu]}(\lambda)=\prod_{m=0}^\infty\left( 
\prod_{j=1}^m\frac{\Gamma(1-\frac1\kappa((\lambda,\alpha)-\frac12(\nu+j\alpha,\alpha)))}
{\Gamma(1-\frac1\kappa((\lambda,\alpha)+\frac12(\nu+j\alpha,\alpha)))}\right)^{d_m^\alpha},$$
where $\Gamma$ is the standard gamma function.

Let $V=V_1\otimes\cdots\otimes V_n$ be a tensor product of finite dimensional $sl_N$ modules.
Set $\Lambda_i(\lambda)=\tr_{V[\nu]}\lambda^{(i)}$, $\varepsilon_{i,j}=\tr_{V[\nu]}\Omega^{(i,j)}$,
$\gamma_i=\sum_{j,\,j\ne i}\varepsilon_{i,j}$, where $i,j=1,\ldots,n$. Set
$$ D_{V[\nu]}(z_1,\ldots,z_n,\lambda)=\prod_{i=1}^nz_i^{\frac1\kappa(\Lambda_i(\lambda)-\frac12\gamma_i)}
\prod_{1\leq i< j\leq n}(z_i-z_j)^{\frac{\varepsilon_{i,j}}{\kappa}}
\prod_{\alpha\in\Sigma_+}X_{\alpha,V[\nu]}(\lambda).$$

Fix a basis $v_1,\ldots, v_d$ in a weight subspace $V[\nu]$. Suppose that $u_i(z,\lambda)=\sum_{j=1}^d
u_{i,j}v_j$, $i=1,\ldots,d$, is a set of $V[\nu]$ valued solutions of the trigonometric KZ
equations and the dynamical equations.  
\begin{theorem}[Corollary 19 in\cite{TV}] \label{t-detTV}
$$\det (u_{i,j})_{1\leq i,j\leq d}=C_{V[\nu]}(\lambda)D_{V[\nu]}(z,\lambda),$$ where $C_{V[\nu]}(\lambda)$
is a function of $\lambda$ (depending also on $V_1,\ldots,V_n$, $\nu$ and $\kappa$) such that
$C_{V[\nu]}(\lambda)$ is $ P^\vee$-periodic, 
$C_{V[\nu]}(\lambda+\kappa\omega^\vee)=C_{V[\nu]}(\lambda)$ for all $\omega^\vee\in P^\vee$.
\end{theorem}

{\bf Remark.} Theorem~\ref{main_sln} implies that Theorem~\ref{t-detTV} can be applied to any set, $u_1,\ldots,u_d$, of
hypergeometric solutions of the trigonometric KZ equations, thus giving a formula for the determinant
of hypergeometric integrals.

{\bf Continuation of the example.} For  $\gg=sl_2$, $n=1$, $V_1=L_p$, we have
$$X_{\alpha,L_p}(\lambda)=\prod_{j=1}^m\frac{\Gamma(1-\frac1\kappa((\lambda,\alpha)-\frac{p}2+m-j))}
{\Gamma(1-\frac1\kappa((\lambda,\alpha)+\frac{p}2-m+j))},
\qquad z_1^{\frac{\Lambda_1(\lambda)}{\kappa}}=z_1^{\frac{(p-2m)(\lambda,\alpha)}{2\kappa}}.$$
The determinant equals the  complex function given as the expression preceding $e_{2,1}^mv_p$
in formula (\ref{ex_sl2}). Denote this function $u_{1,1}$. The explicit formula for the Selberg integral 
in terms of gamma functions implies
$$u_{1,1}=\frac{z_1^{\frac{(p-2m)(\lambda,\alpha)}{2\kappa}}}{m!}
\prod_{j=0}^{m-1}\frac{\Gamma(1+\frac1\kappa(1+j))
\Gamma(1-\frac1\kappa((\lambda,\alpha)-\frac{p}2+m-j-1)\Gamma(\frac{j-p}{\kappa})}
{\Gamma(1+\frac1\kappa)\Gamma(1-\frac1\kappa((\lambda,\alpha)+\frac{p}2-j)}.
$$
Therefore,
$C_{L_p[(p-2m)\frac{\alpha}2]}(\lambda)=(m!)^{-1}
\prod_{j=0}^{m-1}\Gamma(1+\frac1\kappa(1+j))\Gamma(\frac{j-p}{\kappa})(\Gamma(1+\frac1\kappa))^{-1}$.

\section{$N-1$ special normal orders on the set of positive roots of $sl_N$.}\label{spec_ord}

 Matsuo's construction of hypergeometric solutions  uses a PBW-basis of $U(sl_N)$ 
corresponding to the standard 
order of positive roots. We rewrite those solutions using PBW-bases
of $U(sl_N)$ corresponding to some new $N-1$ special (normal) orders of positive roots.
The special orders are defined below.

\subsection{Special normal orders.}
  A linear order on the set of positive roots $\Sigma_+$ of a simple Lie algebra $\gg$
is called normal if for every triple of positive roots  
$\alpha, \alpha+\beta,\beta\in\Sigma_+$ we have either
$\alpha\succ\alpha+\beta\succ\beta$,  or $\beta\succ\alpha+\beta\succ\alpha$. The 
standard order on the set of positive roots of $sl_N$ is a normal order. 

  Consider a linear order on $\Sigma_+$. A permutation $\sigma$ of $\Sigma_+$ is called
{\it an elementary transformation} if $\sigma$ is a reversal of a certain sub-system
$\Sigma_2\subset\Sigma_+$, where $\Sigma_2$ has rank 2 and all elements of $\Sigma_2$ are 
located side by side in the system $\Sigma_+$. An elementary transformation produces a new linear order.
It is known that every elementary transformation converts one normal order into another, and that
any two normal orders on $\Sigma_+$ can be transformed one into the other by a composition of elementary 
transformations, see \cite{AST} and \cite{Z}.

{\bf Example.} For a root system of type $A_m$ ($sl_{m+1}$) we have only two types of elementary
transformations.
\begin{align}
A_1\oplus A_1: & \cdots,\alpha,\beta,\cdots  \rightarrow  \cdots,\beta,\alpha,\cdots&
\mbox{ if } \alpha+\beta\not\in\Sigma_+;\notag\\
A_2: & \cdots,\alpha,\alpha+\beta,\beta,\cdots  \rightarrow  
\cdots,\beta,\alpha+\beta,\alpha,\cdots &\mbox{ if } \alpha+\beta\in\Sigma_+.\notag
\end{align}

  We compute the change of a PBW-basis under a reversal of type $A_1\oplus A_1$ or 
$A_2$.
\begin{lemma} \label{elem_trans}
Let $\gg$ be a simple Lie algebra and let $\alpha$ and $\beta$ be two positive roots.
Choose arbitrary $f_{\alpha}\in\gg_{\alpha}$ and $f_{\beta}\in\gg_{\beta}$. Then, for any triple of
non-negative integers we have the following two identities in $U(\gnm)$.
\begin{align}
A_1\oplus A_1:&  (-1)^{a+b}\frac{f_{\alpha}^a}{a!}\frac{f_{\beta}^b}{b!}=
(-1)^{a+b}\frac{f_{\beta}^b}{b!}\frac{f_{\alpha}^a}{a!};\notag\\
A_2:& (-1)^{a+b+c}\frac{f_{\alpha}^a}{a!}\frac{[f_{\alpha},f_{\beta}]^c}{c!}\frac{f_{\beta}^b}{b!}
=\sum_{p=0}^{\min(a,b)}\bin{c+p}{p}(-1)^{c}(-1)^{a+b+c-p}\frac{f_{\beta}^{b-p}}{(b-p)!}
\frac{[f_{\beta},f_{\alpha}]^{c+p}}{(c+p)!}\frac{f_{\alpha}^{a-p}}{(a-p)!}.\notag
\end{align}
\end{lemma}
{\it Proof.} The proof is a straightforward induction on $a$.\qed \\

Recall that the set $\Sigma_+$ of $sl_N$ is  
$\{\alpha_{k,l}\}_{k<l}$. 
The standard order on $\Sigma_+$ is $\alpha_{k,l}\succ \alpha_{k',l'}$ 
if and only if $l>l'$, or $l=l'$ and $k>k'$. 

\begin{figure}[h]
\begin{center}
\scalebox{.6}{\includegraphics{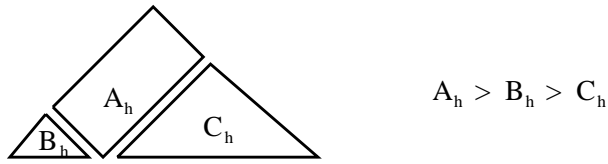}}
\end{center}
\caption{General view of $\succ_h$}\label{Fi:sln_j1}
\end{figure}

\begin{figure}[h]
\begin{center}
\scalebox{.7}{\includegraphics{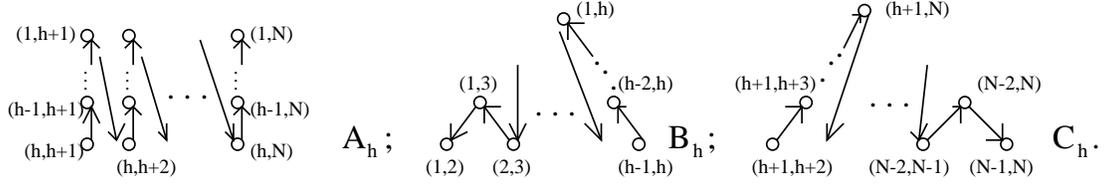}}
\end{center}
\caption{The order  $\succ_h$ within $A_h$, $B_h$, $C_h$. 
Arrows show immediate predecessors.}\label{Fi:sln_j2}
\end{figure}

For $h=1,\ldots,N-1$, define the  
index sets $A_h=\{\alpha_{k,l}\}_{k\leq h<l}$, $B_h=\{\alpha_{k,l}\}_{l\leq h}$,
$C_h=\{\alpha_{k,l}\}_{h<k}$. Define a linear order $\succ_h$ on $\Sigma_+$ by the following rules.
\begin{itemize}
\item $\alpha\succ_h\alpha'\succ_h\alpha''$ for all
 $\alpha\in A_h$, $\alpha'\in B_h$, $\alpha''\in C_h$.
\item The order within the index set $A_h$ is defined by
$\alpha_{k,l}\succ_h \alpha_{k',l'}$ if and only if $l<l'$, or $l=l'$ and $k>k'$ 
for $(k,l), (k',l')\in A_h$.
\item The order within the index set $B_h$ is the standard one.
\item The order within the index set $C_h$ is the opposite to the standard one. 
That is  $\alpha_{k,l}\succ_h \alpha_{k',l'}$ 
if and only if $k<k'$, or $k=k'$ and $l<l'$ for $(k,l), (k',l')\in C_h$.
\end{itemize}
It is easy to see that $\succ_h$ is a normal order,
and  $\succ_{N-1}$ is the standard order. See Figures~\ref{Fi:sln_j1} and 
\ref{Fi:sln_j2} for a pictorial description. 

\subsection{Hypergeometric solutions corresponding to the PBW-basis of type $\succ_h$.}\label{HGsol_h}
Fix $h\in\{1,\ldots,N-1\}$.
Let $V'=V_1\otimes\cdots\otimes V_{n+1}$, where $V_j$ is a highest weight $sl_N$--module with
highest weight $\Lambda_j$ and highest weight vector $v_j$. Fix a weight subspace $V'[\nu']\subset V'$,
$\nu=\sum_{j=1}^n\Lambda_j-\sum_{k=1}^{N-1} m_k\alpha_k$, where 
 $\bm=\sum_{k=1}^{N-1} m_k\alpha_k\in Q^+$. Consider the complex space $\cit^{n+1}$ with coordinates
$z_1,\ldots,z_{n+1}$. 

\paragraph{\bf The PBW--basis of $U(\gnm)$ corresponding to $\succ_h$.}  
For any positive root $\alpha_{k,l}$,  set $a_{k,l}(h)$ 
equal to the number of integers $p$, $k<p< l$, 
such that $\alpha_{k,p}\succ_h\alpha_{p,l}$. The number $a_{k,l}(h)$ counts the $A_2$ subsystems of 
the  order $\succ_h$  which have  $\alpha_{k,l}$ as their middle root 
and which are ordered oppositely to the 
standard order.  Set $F_{k,l}(h)=(-1)^{a_{k,l}(h)}e_{l,k}$. $F(h)=\{-F_{k,l}(h)\}_{k<l}$ is
a basis of $\gnm$ called the basis corresponding to  $\succ_h$.
Order the basis according to $\succ_h$,
$-F_{k,l}(h)\succ -F_{k',l'}(h)$ if and only if $\alpha_{k,l}\succ_h\alpha_{k',l'}$.
The corresponding PBW-basis of $U(\gnm)$ is
\begin{equation}
\left\{ F_{I_0}(h)=(-1)^{\sum i_{l,k}}\frac{F_{h,h+1}(h)^{i_{h+1,h}}}{i_{h+1,h}!}\cdots
\frac{F_{N-1,N}(h)^{i_{N,N-1}}}{i_{N,N-1}!}\right\},
\end{equation}
where $I_0=\{i_{l,k}\}_{k<l}$ runs over all sequences of positive integers.

Let $\{F_{I_1}(h),\ldots,F_{I_{n+1}}(h)\}$ be elements of the PBW-basis of
$U(\gnm)$ corresponding to $\succ_h$.
Set $I=(I_1,\cdots,I_{n+1})$ and $F_I(h)v=F_{I_1}(h)v_1\otimes\cdots\otimes F_{I_{n+1}}(h)v_{n+1}$.
The vector $F_I(h)v$ belongs to $V'[\nu']$ if $I\in P(\bm,n+1)$.

\paragraph{\bf The rational function $\phi(z,t;h)$.} 
For any $I\in P(\bm,n+1)$ and any $s=(j,k,l,q)\in S(I)$, set

\begin{align}
\phi^{(s)}_h&=f^{(s)}\frac{(-1)^{l-1-h}}{t_h^{(\beta_h(I)(s))}-z_j} \,\,\qquad\mbox{ if } k\leq h<l,
\\
\phi^{(s)}_h&=f^{(s)}\frac{1}{t_{l-1}^{(\beta_{l-1}(I)(s))}-z_j} \qquad\mbox{ if } k<l\leq h ,
\notag\\
\phi^{(s)}_h&=f^{(s)}\frac{(-1)^{l-1-k}}{t_k^{(\beta_k(I)(s))}-z_j} \,\,\qquad\mbox{ if } h<k<l,\notag
\end{align}
where $\{\beta_p(I)\}_p$ is the set of bijections we fixed in the definition of $f^{(s)}$, see 
Section~\ref{HGsol}. Set
$$\phi(I,h)=\prod_{s\in S(I)}\phi^{(s)}_h,\qquad \phi(z,t;h)=\sum_{I\in P(\bm,n)}\phi(I,h)F_I(h)v.$$
Notice that $\phi(I,N-1)=\phi(I)$.

\paragraph{\bf Example.}  Let $\gg=sl_3$, $n=0$. The  basis $F(1)=\{-e_{3,2},-(-e_{3,1}),-e_{2,1}\}$ of $\gnm$ 
corresponds to $\succ_1$.
If $ \nu'=\Lambda_1-\alpha_1-\alpha_2$, then
 $$\phi'(z,t;1)=\frac{1}{(t_1^{(1)}-z_1)(t_2^{(1)}-z_1)}e_{2,1}e_{3,2}v_1+
\frac1{(t_2^{(1)}-t_1^{(1)})(t_1^{(1)}-z_1)}e_{3,1}v_1.$$
If $\nu'=\Lambda_1-2\alpha_1-\alpha_2$, then
$$\phi'=\frac{-1}{(t_1^{(1)}-z_1)(t_1^{(2)}-z_1)(t_2^{(1)}-z_1)}\frac{e_{2,1}^2}{2}e_{3,2}v_1+
 \frac{-1}{(t_2^{(1)}-t_1^{(1)})(t_1^{(1)}-z_1)(t_1^{(2)}-z_1)}e_{2,1}e_{3,1}v_1.$$
The corresponding expressions for $\phi'(z,t;2)$ are given in the example in Section~\ref{HGsol},
since $\succ_2$ is the standard order of positive roots of $sl_3$. $\qed$

 The next theorem is our second main result.
\begin{theorem}\label{mat_trig_j}For any $h$, we have 
$\phi(z,t;h)=\phi(z,t)$.
\end{theorem}
The proof of the theorem is in Section~\ref{Pf.5.2}.

  As a corollary, we obtain new  Matsuo's type formulae for solutions of the KZ equations. 
\begin{corollary}
(a)  Let $u(z')$ be the hypergeometric solution of the 
rational KZ equations with parameter $\kappa\in\cit$ in 
the subspace of singular vectors of $V'[\nu']$ indicated in Theorem~\ref{t_mat}.
Then, for every $h=1,\ldots,N-1$, we have   
   $$u(z')=\int_{\gamma(z')}\Phi'(z',t)^{1/k}\left(\sum_{I\in P(\bm,n+1)}\phi(I,h)F_I(h)v'\right)\,dt.$$
(b)  Let $v(z,\lambda)$ be the hypergeometric solution of the trigonometric KZ equation 
with parameters $\kappa\in \cit$ and $\lambda\in\gh$ taking values in the weight space $V[\nu]$ of $V$
indicated in Corollary~\ref{mat_trig}.  For every $h=1,\ldots,N-1$, define a function $v(z,\lambda;h)$ 
$$v(z,\lambda;h):=\int_{\gamma(z)}\Phi(z,t;\lambda)^{1/k}
\left(\sum_{I\in P(\bm,n)}\phi(I,h)F_I(h)v\right)\,dt. $$
Then $v(z,\lambda)=v(z,\lambda;h)$.
\end{corollary}

\section{Additive form of the dynamical difference operators.}\label{add_sect}
\subsection{Statement of the result.} 
Consider a PBW-basis $F=\{F_{I_0}\}_{I_0}$ of $U(\gnm)$. 
Let $\lambda\in\gh$ be generic, and let $M_\lambda$ be the highest weight Verma module with 
highest weight $\lambda$ and highest weight vector $v_\lambda$. The Shapovalov form induces an isomorphism
$S_\lambda:M_\lambda\rightarrow M_\lambda^*$. The set $\{F_{I_0}v_\lambda\}_{I_0}$ is a basis of $M_\lambda$.
Let $\{(F_{I_0}v_\lambda)^*\}_{I_0}$ be the dual basis of $M_\lambda^*$. For every $I_0$, 
there exists a unique element $P_{I_0}(F,\lambda)\in U(\gnm)$ such that $ P_{I_0}(F,\lambda))v_{\lambda}=
S^{-1}_\lambda((F_{I_0}v_\lambda)^*)$. By definition, $P_{I_0}(F,\lambda)$ is a (rational) function
$\gh\rightarrow U(\gnm)$, where $\lambda\mapsto P_{I_0}(F,\lambda)$. 
\paragraph{\bf Remark.} $P_{I_0}(F,\lambda)$ is characterized by the property,
$\tau(P_{I_0}(F,\lambda))v_\lambda^*=(F_{I_0}v_\lambda)^*$.

The map $S^{-1}_\lambda$ corresponds to an element of $M_\lambda\hat{\otimes} M_\lambda$. 
Identify $S^{-1}_\lambda$ with this element. In terms of the basis 
$F$ we have $S^{-1}_\lambda=\sum_{I_0} F_{I_0}v_\lambda\otimes P_{I_0}(F,\lambda)v_\lambda$.
\paragraph{\bf Example.} Let $\gg=sl_2$. Then, $F=\{e_{2,1}^k\}_k$ is a basis of $U(\gnm)$, and
\begin{equation}\label{e-p_2,1^k}
P_{e_{2,1}^k}(F,\lambda)=\frac{e_{2,1}^k}{k!\lambda(\lambda-1)\ldots(\lambda-k+1)}.
\end{equation}
Explicit formula for the element $S^{-1}_\lambda$ for $sl_3$ is given in Appendix C.

Let $V=V_1\otimes\cdots\otimes V_n$ be a tensor product of highest weight $sl_N$ modules and
$V[\nu]\subset V$ a weight subspace. Recall that $\tau:U(\gg)\rightarrow U(\gg)$ is the Chevalley involution, and
$A:U(\gg)\rightarrow U(\gg)$ is the antipode map.

The next theorem is our third main result.
\begin{theorem} \label{main-add}
For every $r=1,\ldots,N-1$, and every $v\in V[\nu]$,  we have
\begin{equation}\label{eq-main-add} \B_{\omega_{[r]},V}(\lambda+\rho+\frac12\nu)v
=\sum_{I_0\in \mathcal{A}(r)}A(F_{I_0}(r))\tau( P_{I_0}(F(r),\lambda)v,
\end{equation}
where $F(r)=\{F_{I_0}(r)\}_{I_0}$ is the PBW-basis of $U(\gnm)$ corresponding to the order $\succ_r$,
and $\mathcal{A}(r)=\{I_0=\{i_{l,k}\}_{k<l}\, |\, i_{l,k}=0 \mbox{ if } k> r, \mbox{ or } l\leq r\}$.
\end{theorem}
Notice that $A(F_{I_0}(r))\in U(\gnm)$ and $\tau( P_{I_0}(F(r),\lambda)\in U(\gnp)$, while according to
the definition, given in Section~\ref{sec2}, $\B_{\omega_{[r]},V}$ is a product of elements $\B^\alpha_V$
each of which contains $U(\gnm)$ and $U(\gnp)$ terms. See Proposition~\ref{BTV-pres} as well.

The theorem is proved in Section~\ref{pf-main-add}.

{\bf Example.} (a)  Let $\gg=sl_2$. Then $\B_{\omega[1],V}(\lambda)=\B^{\alpha_1}(\lambda)$.

(b) Let $\gg=sl_3$.  By definition, we have
$\B_{\omega[1],V}(\lambda)=\B^{\alpha_1+\alpha_2}(\lambda)\B^{\alpha_1}(\lambda)$ and
$\B_{\omega[2],V}(\lambda)=\B^{\alpha_1+\alpha_2}(\lambda)\B^{\alpha_2}(\lambda)$.
Set $(\lambda,\alpha_j)=\lambda_j$, for $j=1,2$, and
$p_k(t)=t(t-1)\ldots(t-k+1)\in\cit[t]$ for any $k\in\nit$. Then, for any $v\in V[\nu]$, we have
\begin{align} 
\B_{\omega[2],V}(\lambda+\rho+\frac12\nu)v&=
\sum_{s,k=0}^{\infty}\frac{e_{3,1}^ke_{3,2}^s}{k!s!}
(\sum_{j=0}^k \frac{ e_{1,2}^{j}e_{1,3}^{k-j}e_{2,3}^{s+j}\bin{k}{j}}{p_{s+j}(\lambda_2)
p_k(\lambda_2+\lambda_1+1)})v,\notag\\
\B_{\omega[1],V}(\lambda+\rho+\frac12\nu)v&=
\sum_{s,k=0}^{\infty}\frac{e_{3,1}^ke_{2,1}^s}{k!s!}
(\sum_{j=0}^k \frac{(-1)^j e_{2,3}^{j}e_{1,3}^{k-j}e_{1,2}^{s+j}\bin{k}{j}}{p_{s+j}(\lambda_1)
p_k(\lambda_2+\lambda_1+1)})v. \notag
\end{align}



\subsection{The universal fusion matrix.} Let $R(\gh)$ be the quotient field of the algebra $U(\gh)$. Then $R(\gh)$
is isomorphic to the field of rational functions over $\gh$, using the identification of the Cartan subalgebra 
$\gh$ with its dual. Define algebras 
$U_1(\gh\oplus\gn_\pm)=U(\gh\oplus\gn_\pm)\otimes_{U(\gh)}R(\gh)$. 
Consider the space $U_1(\gh\oplus\gnm)\hat{\otimes} U_1(\gh\oplus\gnp)$, which is a formal series
completion of $U_1(\gh\oplus\gnm)\otimes U_1(\gh\oplus\gnp)$. For any basis $\{a_K\otimes b_K\}_K$ of
$U(\gnm)\otimes U(\gnp)$,  an element of $U_1(\gh\oplus\gnm)\hat{\otimes} U_1(\gh\oplus\gnp)$
has the form $\sum_{K} a_K\psi_{1,K}\otimes b_K\psi_{2,K}$, where 
$\psi_{1,K},\psi_{2,K}\in R(\gh)$.

Denote $\pi_\gh$ the natural projection $U_1(\gh\oplus\gnm)\hat{\otimes}U_1(\gh\oplus\gnp)\rightarrow
R(\gh)\otimes R(\gh)$. Let $\{x_k\}$ be an orthonormal basis of $\gh$.
\begin{theorem}[Proposition 1 in \cite{ABBR}]\label{t-ABBR}
 For every $\lambda\in\gh$,
there exists a unique solution $J(\lambda)$ of the ABBR-equation
\begin{equation}\label{ABBReq}
[1\otimes(\lambda+\rho-\frac12\sum_k x_k^2),J(\lambda)]=
-(\sum_{\alpha\in\Sigma_+}e_{-\alpha}\otimes e_\alpha)J(\lambda),
\end{equation}
such that $J(\lambda)$ belongs to $U_1(\gh\oplus\gnm)\hat{\otimes} U_1(\gh\oplus\gnp)$ and 
$\pi_\gh(J(\lambda))=1\otimes 1$. Moreover, $F(\lambda)$ is of weight zero, and has the following expansion
$$ J(\lambda)=\sum_{K=0}^\infty a_K\otimes(b_K\psi_K(\lambda)),$$
where $a_K\in U(\gnm)$, $b_K\in U(\gnp)$, $\psi_K\in R(\gh)$, for all K, $a_0=b_0=\psi_0(\lambda)=1$.
\end{theorem}
For technical reasons a proof of the theorem is given at the end of this section.

The solution $J(\lambda)$ is called the universal  fusion matrix of $U(\gg)$. Theorem~\ref{c-add-inter} and
equality~(\ref{e-p_2,1^k}) give a formula 
for the universal  fusion matrix of $U(sl_2)$. See Appendix C for a similar 
formula in the $sl_3$ case.

Set $Q^\dagger(\lambda)=\sum_{K\geq 0} A(a_K)b_K\psi_K(\lambda)$, cf. \cite{EV}. The 
 element $Q^\dagger(\lambda)$ belongs to a formal series completion of 
$U(sl_N)\otimes_{U(\gh)}R(\gh)$.
Its action is well defined in a tensor product of highest weight $sl_N$ modules.

\begin{theorem}[Theorem 34 in \cite{EV}]  \label{EV-T34}
Let $w_0$ be the longest element in the Weyl group $\W$. For any $v\in V[\nu]$, we have 
$$ \B_{w_0,V}(\lambda+\rho-\frac12\nu)v=Q^\dagger(\lambda)v.$$
\end{theorem}

Next, we present a connection between the Shapovalov form and the universal  fusion matrix.
\begin{theorem}[Prop. 13.1 in \cite{ES}] \label{prop.13.1}
We have $S^{-1}_\lambda=\left((1\otimes\tau)J(0)\right)(v_\lambda\otimes v_\lambda)$,
where $S^{-1}_\lambda$ is identified with the corresponding element of $M_\lambda\hat{\otimes} M_\lambda$.
\end{theorem}
More explicitly, in terms of a basis $F$ of $U(\gnm)$, we have
\begin{equation}\label{J(0)}
\sum_{I_0} F_{I_0}v_\lambda\otimes P_{I_0}(F,\lambda)v_\lambda=
\left(\sum_{K\geq 0} a_K\otimes \tau(b_K\psi_K(0))\right)(v_\lambda\otimes v_\lambda).
\end{equation}  

 Let $F(\lambda)$ be a rational function of $\lambda\in\gh$ with values in $R(\gh)$.
Assume that the $\lambda$ dependence in $F(\lambda)$ is through a finite number of inner products 
$\{(\lambda,\beta_j)\}_{j=1}^k$, $F(\lambda)=F((\lambda,\beta_1),\ldots,(\lambda,\beta_k))$, where
$\beta_1,\ldots,\beta_k\in\gh$. Set $F_+(\lambda)=F((\lambda,\beta_1)+\beta_1,\ldots,(\lambda,\beta_k)+\beta_k)$. 
The function $F_+(\lambda)$ is a rational function of $\lambda\in\gh$ with values in 
$R(\gh)$.
\begin{theorem} \label{c-add-inter}
Let $F=\{F_{I_0}\}_{I_0}$ be a homogeneous basis of $U(\gnm)$.  
Set $J_+(\lambda)=\sum_K a_K\otimes (b_K(\psi_K)_+(\lambda))$.
Then $J_+(\lambda)$ is well defined, and
$$J_+(\lambda)=\sum_{I_0} F_{I_0} \otimes \tau(P_{I_0}(F,\lambda)).$$
\end{theorem}
Notice that, in particular, the theorem says that $J_+(\lambda)$ takes values in 
$U(\gnm)\hat{\otimes} U(\gnp)$ for any $\lambda\in\gh$.
The proof of the theorem is given at the end of this section.

 We derive the following corollary from Theorems~\ref{EV-T34} and \ref{c-add-inter}.
\begin{corollary} \label{add_w0}
Let $w_0$ be the longest element in the Weyl group $\W$. 
Then,  $$\B_{w_0,V}(\lambda+\rho+\frac12 \nu)v=\left( \sum_{I_0} A(F_{I_0})\tau(P_{I_0}(F,\lambda))\right)v,$$
for any homogeneous basis $F=\{F_{I_0}\}_{I_0}$ of $U(\gnm)$, and any  homogeneous element $v\in V[\nu]$.
\end{corollary}
Indeed, by Theorem~\ref{EV-T34}, we have
$$\B_{w_0,V}(\lambda+\rho+\frac12\nu)v=Q^\dagger(\lambda+\nu)v=
\sum_{K\geq 0} A(a_K)b_K\psi_K(\lambda+\nu)v.$$
Theorem~\ref{c-add-inter} gives
$$J_+(\lambda)=
\sum_{K\geq 0} a_K\otimes (b_K(\psi_K)_+(\lambda))=\sum_{I_0} F_{I_0} \otimes \tau(P_{I_0}(F,\lambda)).$$
Finally, the equality $((\lambda,\beta)+\beta)v=((\lambda,\beta)+(\nu,\beta))v=(\lambda+\nu,\beta)v$ for any
$\lambda,\beta\in \gh$, implies
$$\sum_{K\geq 0} A(a_K)b_K\psi_K(\lambda+\nu)v=\sum_{K\geq 0} A(a_K)b_K(\psi_K)_+(\lambda)v
=\sum_{I_0}A(F_{I_0})\tau(P_{I_0}(F,\lambda))v.\qquad\qed$$

{\it Proof of Theorem~\ref{t-ABBR}.} Fix $\lambda\in \gh$. Denote $wt(\,.\,)$ 
the weight function defined on homogeneous elements in $U(\gg)$ with values in $\gh$. 
Fix a basis $\{a_K\otimes b_K\}_K$ of $U(\gnm)\otimes U(\gnp)$, such that $a_K$ 
and $b_K$, $K\geq 0$,  are homogeneous elements of $U(\gnm)$ and $U(\gnp)$, respectively.
  Let $J(\lambda)=\sum_{K\geq 0}a_K\psi_{1,K}(\lambda)\otimes b_K\psi_{2,K}(\lambda)$, 
where $\psi_{1,K}(\lambda),\psi_{2,K}(\lambda)\in R(\gh)$

The ABBR equation~(\ref{ABBReq}) is  
\begin{align}
\sum_{K\geq 0} (a_K \psi_{1,K}(\lambda)&\otimes b_K \psi_{2,K}(\lambda) )
\left( 1\otimes((\lambda+\rho,wt(b_K))-\frac12(wt(b_K),wt(b_K))-wt(b_K))\right)=\notag\\
 &=-(\sum_{\alpha\in\Sigma_+} \sum_{K\geq 0} (e_{-\alpha}a_K \psi_{1,K}(\lambda)
\otimes e_\alpha b_K \psi_{2,K}(\lambda)).\label{addB_i1}
\end{align}
Equality~(\ref{addB_i1}) presents a system of recurrence relations for  
$\{\psi_{1,K}(\lambda)\otimes \psi_{2,K}(\lambda)\}_{k\geq 0}\subset R(\gh)\otimes R(\gh)$. 
Together with the initial conditions $\psi_{1,0}(\lambda)=\psi_{2,0}(\lambda)=1$, 
they uniquely determine $\psi_{1,K}(\lambda)\otimes \psi_{2,K}(\lambda)$ for all $K$. Moreover,
the solution is such that $\psi_{1,K}(\lambda)=1$. Set $\psi_K(\lambda)=\psi_{2,K}(\lambda)$.
$\qed$

{\it Proof of Theorem~\ref{c-add-inter}.} Observe that $\psi_K(\lambda)$ is a function of $\lambda$ 
with values in $R(\gh)$, it  depends on $\lambda$ via the inner products 
$\{(\lambda,wt(b_L(\lambda))\}$, $0\leq wt(b_L(\lambda))\leq wt(b_K(\lambda))$.
Thus, the shifted function $(\psi_K)_+(\lambda)$ is well 
defined for any $K$, and  $J_+(\lambda)$ is well defined.

Write a system of recurrence relations for $J_+(\lambda)$ by shifting the system of reccurence relations
(\ref{addB_i1}). We have
\begin{align}
\sum_{K\geq 0} (a_K &\otimes (b_K(\psi_K)_+(\lambda)) )
\left( 1\otimes((\lambda+\rho,wt(b_K))-\frac12(wt(b_K),wt(b_K)))\right)=\notag\\
\label{addB_i} &=-\sum_{\alpha\in\Sigma_+} \sum_{K\geq 0} 
(e_{-\alpha}a_K \otimes (e_\alpha b_K(\psi_K)_+(\lambda)) ).
\end{align}
Equation~(\ref{addB_i}) implies a system of reccurence relations for the functions 
$\{(\psi_K)_+(\lambda)\}_{K\geq 0}$. The initial condition is $(\psi_0)_+(\lambda)=1$. 
Therefore, for the unique solution, each $(\psi_K)_+(\lambda)$ is a complex-valued functions.
The shifted universal fusion matrix $J_+(\lambda)$ takes values in $U(\gnm)\hat{\otimes} U(\gnp)$. 

Next we study $J(0)$ which is connected to $\sum_{I_0} F_{I_0} \otimes \tau(P_{I_0}(F,\lambda))$ by 
Theorem~\ref{prop.13.1}. The ABBR reccurence relations ~(\ref{addB_i1}) applied to $J(0)$ give
\begin{align}\label{addB_i0}
\sum_{K\geq 0} (a_K &\otimes (b_K\psi_K(0)) )
\left( 1\otimes((\rho,wt(b_K))-\frac12(wt(b_K),wt(b_K))-wt(b_K))\right)=\\
 &=-\sum_{\alpha\in\Sigma_+} \sum_{K\geq 0} (e_{-\alpha}a_K \otimes (e_\alpha b_K\psi_K(0)) ).\notag
\end{align}
Consider $(1\otimes \tau)(J(0))=\sum_{K\geq 0}a_K\otimes \tau(b_K\psi_K(0))$ 
$\in U(\gnm)\hat{\otimes} U_1(\gh\oplus\gnm)$. Its action on the vector 
$v_\lambda\otimes v_\lambda\in M_\lambda\otimes M_\lambda$ defines an element of $M_\lambda\hat{\otimes} M_\lambda$.
Notice that $\tau((\rho,wt(b_K))-\frac12(wt(b_K),wt(b_K))-wt(b_K))=
((\rho,wt(b_K))-\frac12(wt(b_K),wt(b_K))+wt(b_K))$,
because $\tau$ restricted to $\cit$ is the identity and $\tau$ restricted to $\gh$ is multiplication by $-1$.
Since $wt(b_K) v_\lambda=(\lambda,wt(b_K))v_\lambda$, equality~(\ref{addB_i0}) implies
\begin{align}\label{abbr_0}
\sum_{K\geq 0} (a_K v_\lambda &\otimes \tau(b_K\psi_K(0))v_\lambda)((\lambda+\rho,wt(b_K))
-\frac12(wt(b_K),wt(b_K)))=\\
 &=-\sum_{\alpha\in\Sigma_+} \sum_{K\geq 0} e_{-\alpha}a_K v_\lambda\otimes 
\tau(e_\alpha b_K\psi_K(0))v_\lambda.\notag
\end{align}
For any $K$, define a complex-valued function $\psi_K(0,\lambda)$, such that $\tau(\psi_K(0))v_\lambda=
\psi_K(0,\lambda)v_\lambda$.
Therefore,  we have
$ \psi_K(0,\lambda)(a_K v_\lambda\otimes\tau(b_K)v_\lambda)=a_K v_\lambda\otimes\tau(b_K\psi_K(0))v_\lambda$.
Now (\ref{abbr_0}) takes the form
\begin{align}\label{abbr_01}
\sum_{K\geq 0} (a_K v_\lambda &\otimes \tau(b_K)v_\lambda)\psi_K(0,\lambda)((\lambda+\rho,wt(b_K))
-\frac12(wt(b_K),wt(b_K)))=\\
 &=-\sum_{\alpha\in\Sigma_+} \sum_{K\geq 0} (e_{-\alpha}a_K v_\lambda\otimes 
\tau(e_\alpha b_K)v_\lambda)\psi_K(0,\lambda).\notag
\end{align}

Equality~(\ref{abbr_01}) together with the isomorphism between $M_\lambda$ and $U(\gnm)$ imply a
system of recurrence relations for the functions $\{\psi_K(0,\lambda)\}_{K\geq 0}$ with 
initial condition $\psi_0(0,\lambda)=1$. This systems of recurrence relations coincides
with the systems of recurrence relations for $\{(\psi_K)_+(\lambda)\}_{k\geq 0}$. 
The initial condition in both cases is the same. Therefore the solutions coincide, i.e. for every $K$,
$(\psi_K)_+(\lambda)=\psi_K(0,\lambda)$ and $a_K\otimes (b_K(\psi_K)_+(\lambda))=a_K\otimes (b_K\psi_K(0,\lambda))$.

In our notation, Theorem~\ref{prop.13.1} together with the linear isomorphism  between $M_\lambda$ and $U(\gnm)$ 
give the following explicit formula, $\sum_{K\geq 0} a_K \otimes \tau(b_K)\psi_K(0,\lambda)
=\sum_{I_0} F_{I_0} \otimes P_{I_0}(F,\lambda)$. Therefore,
$J_+(\lambda)=\sum_{I_0} F_{I_0} \otimes \tau(P_{I_0}(F,\lambda))$, which
is the statement of the corollary. \qed

\subsection{Normal orders on the set of positive roots and multiplicative presentations of 
$\B_{w_0,V}(\lambda)$.} 
Let $\omega=s_{i_m}\cdots s_{i_1}$ be a reduced decomposition of an element $\omega\in\W$.
Set $\alpha^1=\alpha_{i_1}$,\hspace{.5cm}$\alpha^p=s_{i_1}\cdots s_{i_p-1}\alpha_{i_p}$, for $p=2,\ldots, m$. 
\begin{proposition}[\cite{TV}]\label{BTV-pres}
$\B_{\omega,V}(\lambda)=\B_V^{\alpha^m}(\lambda)\cdots\B_V^{\alpha^1}(\lambda)$.
\end{proposition}

Let $\omega_0=s_{i_m}\cdots s_{i_1}$ be a reduced decomposition of the longest element 
$\omega_0\in\W$. Then $\alpha^m\succ\cdots\succ\alpha_1$
is a normal order on $\Sigma_+$, and any normal order corresponds to a reduced decomposition of the longest 
element, see \cite{Z}.
\begin{corollary} \label{B_0pres}
If $\alpha^m\succ\cdots\succ\alpha^1$ is a normal order on $\Sigma_+$, then
$$\B_{\omega_0,V}(\lambda)=\B_V^{\alpha^m}(\lambda)\cdots\B_V^{\alpha^1}(\lambda).$$
\end{corollary}

\subsection{The proof of Theorem~\ref{main-add}.}\label{pf-main-add}
  Fix $r\in \{1,\ldots,N-1\}$. Recall that 
the element $\omega_{[r]}$ is defined
by $\omega_{[r]}=\omega_0\omega_0^r\in \W$, where $\omega_0$ (respectively, $\omega_0^r$) 
is the longest element in $\W$ (respectively, in $\W^r$ generated by all simple reflections 
$s_l$ preserving  the dual fundamental weight $\omega_r^\vee$). The explicit form of 
$\omega_r^\vee=\sum_{k=1}^r(1-\frac{r}{N})e_{k,k}
-\sum_{k=r+1}^N\frac{r}{N}e_{k,k}$ implies that $\W^r$ is generated by the simple reflections
$s_1,\ldots,s_{r-1},s_{r+1},\ldots,s_{N-1}$. 
Denote $\W_2$ the subgroup of $\W$ generated by
$s_1,\ldots,s_{r-1}$, and denote $\W_1$ the subgroup of $\W$ generated by $s_{r+1},\ldots,s_{N-1}$.
Denote $R_2$ the  root system with base $\alpha_1,\ldots,\alpha_{r-1}$, and
denote $R_1$  the  root system with base $\alpha_{r+1},\ldots,\alpha_{N-1}$. 
Then $\W_1$ is the Weyl group of $R_1$, and
$\W_2$ is the Weyl group of $R_2$, and $\W^h=\W_2\times\W_1$.

  Consider the normal order $\succ_r$. Set $m_1=(N-r-1)(N-r)/2$, $m_2=(r-1)r/2$, $m_3=r(N-r)$, 
$m=m_1+m_2+m_3$. Write explicitly the order $\succ_r$ on $\Sigma_+$,
$\alpha^m\succ_r\cdots\succ_r\alpha^1$.
We have $\Sigma_+=A_r\cup B_r\cup C_r$. Moreover, $A_r=\{\alpha^m,\ldots,\alpha^{1+m_2+m_1}\}$,
$B_r$ is the set of positive roots for $R_2$,
$C_r$ is the set of positive roots for $R_1$, 
and $A_r\succ_r B_r\succ_r C_r$. 
 Let $\omega_0=s_{i_{m}}\cdots s_{i_1}$ be the reduced decomposition of $\omega_0$ 
corresponding to $\succ_r$.
Denote $\omega_1=s_{i_{m_1}}\cdots s_{i_1}$, $\omega_2=s_{i_{m_2+m_1}}\cdots s_{i_{1+m_1}}$,
$\omega_3=s_{i_{m}}\cdots s_{i_{1+m_2+m_1}}$. Clearly, $\omega_2\omega_1=\omega_0^r$, 
$\omega_3=\omega_{[r]}$.

  Let $\{\beta^{m_3},\ldots,\beta^1\}$ be the set of positive roots 
corresponding to the reduced presentation $\omega_{[r]}=s_{i_{m}}\cdots s_{i_{1+m_2+m_1}}$.
Then $\B_{\omega_{[r]},V}(\lambda)=\B_V^{\beta^{m_3}}(\lambda)\cdots\B_V^{\beta^1}(\lambda)$. 
The definition of the correspondence between normal orders and reduced presentations of $\omega_0$
implies that $A_r=$ $(\omega_0^r)^{-1}(\{\beta^{m_3},\ldots,\beta^1\})$ as ordered sets. Thus,
we have $\{\beta^{m_3},\ldots,\beta^1\}=$ \\$\omega_0^r(\{\alpha^m,\ldots,\alpha^{1+m_2+m_1}\})$. 
\begin{lemma} \label{inverse_A_h}
The reflection $\omega_0^r$ reverses $A_r$,
  $\omega_0^r(A_r)=\{\alpha^{1+m_1+m_2},\ldots,\alpha^m\}$. Therefore,
$$ \B_{\omega_{[r]},V}(\lambda)=\B_V^{\alpha^{1+m_1+m_2}}(\lambda)\cdots\B_V^{\alpha^m}(\lambda).$$
\end{lemma}
The proof is a straightforward computation. $\qed$
\begin{corollary}\label{add-prod}
\begin{equation*}
\B_{\omega_0,V}(\lambda)=\B_{\omega_1,V}(\lambda)\B_{\omega_2,V}(\lambda)\B_{\omega_{[r]},V}(\lambda).
\end{equation*}
\end{corollary}
{\it Proof.}
Consider the normal order $\prec^r$ on $\Sigma_+$ reverse to $\succ_r$. For any
$\alpha,\beta\in\Sigma_+$ we have $\alpha\prec^r\beta$ if and only if
$\alpha\succ_r\beta$. Notice that $A_r\prec^h B_r\prec^r C_r$ and the order $\prec^r$ within each of the sets 
$A_r$, $B_r$, $C_r$ is reverse to $\succ_r$. Consider the presentation of
$\B_{\omega_0,V}(\lambda)$ corresponding to $\prec^r$ given in Corollary~\ref{B_0pres}. 
Apply Lemma~\ref{inverse_A_h} to obtain the statement of the corollary. $\qed$

Recall that $F(r)=\{F_{I_0}(r)\}_{I_0}$ is the PBW-basis of $U(\gnm)$ corresponding to $\succ_r$.
\begin{lemma}\label{l-Binter} For any homogeneous vector $v\in V[\nu]$, we have
$$\B_{\omega_{[r]},V}(\lambda+\rho+\frac12 \nu)v
=\sum_{I_0\in \mathcal{A}(r)}A(F_{I_0}(r)) Q_{I_0}(F(r),\lambda)v,$$
where $\{Q_{I_0}(F(r),\lambda)\}_{I_0}$ is a subset of homogeneous elements of $U(\gnp)$.
\end{lemma}
{\it Proof.} For every positive root $\alpha$, every $\lambda\in \gh$ and every positive integer $s$, set
$p_{0,\lambda,\alpha}=1$, and
$p_{s,\lambda,\alpha}=s!((\lambda+\rho,\alpha)-1)\ldots((\lambda+\rho,\alpha)-s)\in \cit$. Use formula 
(\ref{p-oper}) to obtain 
$$\B_V^{\alpha}(\lambda+\rho+\frac12 \nu)v=\sum_{s=0}^\infty 
F_\alpha^sE_\alpha^s (p_{s,\lambda,\alpha})^{-1}v.$$
 Since by definition 
$A_r=\{\alpha_{r,r+1},\alpha_{r-1,r+1},\ldots,\alpha_{2,N},\alpha_{1,N}\}$, see Figure~\ref{Fi:sln_j2}, 
Lemma~\ref{inverse_A_h} implies
\begin{equation}\label{B-inter}
\B_{\omega_{[r]},V}(\lambda+\rho+\frac12 \nu)v=\sum_{s_{r,r+1}=0}^\infty\ldots
\sum_{s_{1,N}=0}^\infty \frac{F_{\alpha_{1,N}}^{s_{1,N}}E_{\alpha_{1,N}}^{k_{1,N}}\ldots
F_{\alpha_{r,r+1}}^{s_{r,r+1}}E_{\alpha_{r,r+1}}^{s_{r,r+1}}}
{p_{s_{k_{1,N},\lambda,\alpha_{1,N}}}\ldots p_{s_{r,r+1},\lambda,\alpha_{r,r+1}}}v.
\end{equation}
Next, we need to rearrange the right hand side of formula~(\ref{B-inter}) by 
moving all elements of type $F_{k,l}$ 
to the left through elements of type $E_{k',l'}$, using the commutation relations of $sl_N$. 
Since $E_{\alpha_{k',l'}}=e_{k',l'}$ and $F_{\alpha_{k,l}}=e_{l,k}$, we have
\begin{align}
& (a)\quad[E_{\alpha_{k',l}},F_{\alpha_{k,l}}]=E_{\alpha_{k',k}}, \mbox{ if } k'<k; \qquad
(b)\quad[E_{\alpha_{k',l}},F_{\alpha_{k,l}}]= F_{\alpha_{k,k'}}, \mbox{ if } k<k';\notag \\
& (c)\quad[E_{\alpha_{k,l'}},F_{\alpha_{k,l}}]=-E_{\alpha_{l,l'}}, \mbox{ if } l<l'; \qquad
(d)\quad[E_{\alpha_{k,l'}},F_{\alpha_{k,l}}]= - F_{\alpha_{l',l}}, \mbox{ if } l'<l;\notag\\
&(e)\quad [E_{\alpha_{k,l}},F_{\alpha_{k,l}}]=H_{\alpha_{k,l}}; \qquad
(f)\quad[E_{\alpha_{k',l'}},F_{\alpha_{k,l}}]=0,  \mbox{ if } k\ne k' \mbox{ and } l\ne l'.\notag
\end{align}

The initial ordering of the product of $E$'s and $F$'s in formula (\ref{B-inter}) is according
to the $\succ_r$ ordering of $A_r$, $F_{\alpha_{k,l}}^{s_{k,l}}E_{\alpha_{k,l}}^{s_{k,l}}$ is to the 
left of $F_{\alpha_{k',l'}}^{s_{k',l'}}E_{\alpha_{k',l'}}^{s_{k',l'}}$ if and only if 
$\alpha_{k',l'}\succ_r\alpha_{k,l}$. See Figure~\ref{Fi:sln_j2} for the definition of $\succ_r$ on $A_r$.
 The set of indices $\{(k,l)\}$ is characterized by  $1\leq k\leq r<l\leq N$. 

We first move all $F_{\alpha_{2,N}}$'s to the left through all $E_{\alpha_{1,N}}$'s, 
then we move all $F_{\alpha_{3,N}}$'s and so on. The last step is to move all $F_{\alpha_{r,r+1}}$'s 
to the left through all $E$'s resulting from the previous steps. 

{\bf Claim.} After each step we have:\\
(*) all $F_{\alpha_{k,l}}$'s in the respective rearranged version of formula 
(\ref{B-inter}) satisfy $k\leq r <l$; \\
(**) if $k=k'\leq r <l'\leq l$, or $k\leq k'\leq r<l'=l$ then 
any $E_{\alpha_{k',l'}}$ is to the right of any $F_{\alpha_{k,l}}$.

In the initial state represented in   formula (\ref{B-inter}) the claim follows from the definition of the order
$\succ_r$. Now, assume that the claim is valid after $p$ rearranging steps. Assume that we move
$F_{\alpha_{k,l}}$ at the $p+1$-st step.

If the $p+1$-st step is of type $(a)$ or $(c)$, then we get an additional $\pm E_{\alpha_{k',l'}}$ with
$k'<l'\leq r$ or $r<k'<l'$. $E$'s of this type does not affect the claim, so it remains valid.

If the $p+1$-st step is of type $(b)$ or $(d)$, then we get an additional $\pm F_{\alpha_{k',l'}}$ with
$k\leq k'\leq r<l'\leq l$. Other $F$'s are excluded because of condition (**). Moreover,  a simple check up
shows that condition (**) remains valid for the additional summand.

Condition (**) shows that the $p+1$-st step can  not be of type $(e)$, and
after a step of type $(f)$ the claim is trivially satisfied. 

  Inductive argument shows that the claim is satisfied after each step of the rearrangement. As a final result 
we get 
\begin{equation}\label{B-inter1}
\B_{\omega_{[r]},V}(\lambda+\rho+\frac12 \nu)v=\sum_{s_{r,r+1}=0}^\infty\ldots
\sum_{s_{1,N}=0}^\infty (\prod_{k\leq r<l} F_{\alpha_{k,l}}^{s_{k,l}})\widetilde{Q}_{\{s_{k,l}\}}
(\lambda)v,
\end{equation}
where each $\widetilde{Q}_{\{s_{k,l}\}}(\lambda)$ is a homogeneous element of $U(\gnp)$. 
Finally, identify the set $\{s_{k,l}\in\zit_{\geq 0}\,|\, 1\leq k\leq r <l\leq N\}$ 
with the set  
$\mathcal{A}(r)=\{s_{k,l}\in\zit_{\geq 0}\,|\, 1\leq k<l\leq N,\,s_{k,l}=0 \mbox{ if } r<k, \mbox{ or }  l\leq r\}$.
The relation $(\prod_{k\leq r<l} F_{\alpha_{k,l}}^{s_{k,l}})=\pm (\prod_{k\leq r<l} s_{k,l}!)
A(F_{\{s_{k,l}\}}(r))$ together with equality (\ref{B-inter1}) imply the statement of the lemma. $\qed$

Since $\omega_1$ is the longest element in $\W(N-r-1)$ and $\omega_2$ is the longest element in $\W(r-1)$
we can apply Corollary~\ref{add_w0} to $\B_{\omega_1,V}$ and $\B_{\omega_2,V}$. Choose
$F(r)=\{F_{I_0}(r)\}_{I_0}$ as the weighted PBW-basis of $U(\gnm)$. Then,
for every $v\in V[\nu]$,  we have
\begin{align}
\B_{\omega_1,V}(\lambda+\rho+\frac12 \nu)v &= 
\sum_{I_0\in \mathcal{C}(r)} A(F_{I_0}(r))\tau(P_{I_0}(F(r),\lambda))v,\notag\\
\B_{\omega_2,V}(\lambda+\rho+\frac12 \nu)v &= 
\sum_{I_0\in \mathcal{B}(r)} A(F_{I_0}(r))\tau(P_{I_0}(F(r),\lambda))v,\notag
\end{align} 
where $\mathcal{C}(r)=\{I_0=\{i_{l,k}\}_{k<l}\,|\, i_{l,k}=0 \mbox{ if } k\leq r\}$, and
$\mathcal{B}(r)=\{I_0=\{i_{l,k}\}_{k<l}\,| \, i_{l,k}=0 \mbox{ if } l>r\}$. 
It is easy to see that, for any $I_1\in \mathcal{C}(r)$, $I_2\in\mathcal{B}(r)$, we have 
$[\tau(P_{I_1}(F,\lambda),A(F_{I_2})]=0$. Therefore,
\begin{equation}\label{add_w_0^h}
\B_{\omega_1,V}(\lambda+\rho+\frac12 \nu)\B_{\omega_2,V}(\lambda+\rho+\frac12 \nu)v=
\sum_{I_0\in \mathcal{BC}(r)} A(F_{I_0}(r))\tau(P_{I_0}(F(r),\lambda))v,
\end{equation}
where $\mathcal{BC}(r)=\{I_0=\{i_{l,k}\}_{k<l}\,|\, i_{l,k}=0 \mbox{ if } k\leq r <l\}$.

Apply Corollary~\ref{add_w0} to the right hand side of Corollary~\ref{add-prod}, and 
Lemma~\ref{l-Binter} and formula~(\ref{add_w_0^h}) to the left hand side of Corollary~\ref{add-prod}.
For any $v\in V$, we obtain
\begin{align}
\sum_{I_0} &A(F_{I_0}(r))\tau(P_{I_0}(F(r),\lambda))v=\label{eq-add-prod1}\\
&=\sum_{I_0\in \mathcal{BC}(r)}\hspace{-7pt} A(F_{I_0}(r))\tau(P_{I_0}(F(r),\lambda))
\hspace{-7pt}\sum_{I_0\in \mathcal{A}(r)}\hspace{-7pt}A(F_{I_0}(r)) Q_{I_0}(F(r),\lambda)v.\notag
\end{align}
Set $V=M_\lambda$ the Verma module with highest weight $\lambda$. Weight considerations imply
$$\tau(P_{I_0}(F(r),\lambda))F_{J_0}(r)v_\lambda=0,\quad Q_{I_0}(F(r),\lambda)F_{J_0}(r)v_\lambda=0
\mbox{ if } wt(F_{J_0}(r)v_\lambda)\not \leq wt(F_{I_0}(r)v_\lambda).$$
The defining property of $P_{I_0}(F(r),\lambda)$ gives
$$S_\lambda(\tau(P_{I_0}(F(r),\lambda))F_{J_0}(r)v_\lambda,v_\lambda)=
S_\lambda(F_{J_0}(r)v_\lambda,A(P_{I_0}(F(r),\lambda))v_\lambda)=\pm \delta_{I_0,J_0}.$$  
Therefore $\tau(P_{I_0}(F(r),\lambda))F_{I_0}(r)v_\lambda=
\pm v_\lambda$,  and $\tau(P_{I_0}(F(r),\lambda))F_{I_0}(r)v_\lambda=0$ if $J_0\ne I_0$ and
$wt(F_{J_0}(r)v_\lambda)= wt(F_{I_0}(r)v_\lambda)$.

Finally, apply equality~(\ref{eq-add-prod1}) to $F_{I_0}(r)v_\lambda$, for $I_0\in\mathcal{A}$, and use induction 
on the partialy ordered set $\{ wt(F_{I_0}(r)v_\lambda)\}_{I_0\in\mathcal{A}}$ to obtain
$Q_{I_0}(F(r),\lambda)F_{J_0}(r)v_\lambda=\tau(P_{I_0}(F(r),\lambda))F_{J_0}(r)v_\lambda$, for any $J_0$, and thus
$Q_{I_0}(F(r),\lambda)=\tau(P_{I_0}(F(r),\lambda))$. The last equality and Lemma~\ref{l-Binter} imply the statement 
of Theorem~\ref{main-add}
\qed

\section{Proof of Theorem~\ref{mat_trig_j}.}\label{Pf.5.2}

\subsection{Sequence of elementary transformations converting $\succ_h$ into $\succ_{h-1}$.}
Recall that $h\in\{1,\ldots,N-1\}$. The following construction is very important. Namely,
  there exists a sequence $\sigma$ of elementary transformations 
\begin{equation}\label{ord_j_to_j-1}
\sigma=\{\sigma_1,\ldots,\sigma_M\},
\end{equation}
which converts the normal order  $\succ_h$ into the normal 
order $\succ_{h-1}$. The sequence is such that, 
for every pair of positive integers $(k,l)$, $1\leq k<h<l\leq N$,
the $A_2$ elementary transformation of  type
$\cdots, \alpha_{h,l},\alpha_{k,l},\alpha_{k,h},\cdots\rightarrow 
\cdots, \alpha_{k,h},\alpha_{k,l},\alpha_{h,l},\cdots$   is  used exactly once in this sequence and
all other elementary transformations in the sequence are of type  $A_1\oplus A_1$.
The sequence is given as follows.

\paragraph{\bf The construction of $\sigma$.} We have $\Sigma_+=A_h\cup B_h\cup C_h=
A_{h-1}\cup B_{h-1}\cup C_{h-1}$, and \\$A_h\succ_h B_h\succ_h C_h$, 
\hspace{.5cm}$A_{h-1}\succ_{h-1} B_{h-1}\succ_{h-1} C_{h-1}$. Moreover,
\begin{align}
&A_{h-1}=A_h+ \{\alpha_{k,h}\}_{k=1}^{k=h-1}-\{\alpha_{h,l}\}_{l=h+1}^{l=N},\quad
B_{h-1}=B_h-\{\alpha_{k,h}\}_{k=1}^{k=h-1}, \notag\\
& C_{h-1}=C_h+  \{\alpha_{h,l}\}_{l=h+1}^{l=N}, \notag
\end{align}
where $\{\alpha_{k,h}\}_{k=1}^{k=h-1}\subset B_h$, 
$\{\alpha_{h,l}\}_{l=h+1}^{l=N}\subset A_h$.
The roots $\{\alpha_{k,h}\}_{k=1}^{k=h-1}$ are the largest elements in $B_h$ according to
$\succ_h$ and are the largest elements in $A_{h-1}$ according to $\succ_{h-1}$. 
The linear order $\succ_h$ on the set  $\{\alpha_{k,h}\}_{k=1}^{k=h-1}$ 
is the same as the linear order $\succ_{h-1}$.
The roots $\{\alpha_{h,l}\}_{l=h+1}^{l=N}$ are the largest elements in $C_{h-1}$ according to 
$\succ_{h-1}$. The linear order
$\succ_h$ on the set  $\{\alpha_{h,l}\}_{l=h+1}^{l=N}$ 
is the same as the linear order $\succ_{h-1}$.

The procedure changing the order $\succ_h$ into the order $\succ_{h-1}$ consists
of two parts described in detail below. In part one, we move the sequence of 
roots $\{\alpha_{k,h}\}_{k=1}^{k=h-1}\subset B_h$ to the left 
through the set $A_h$. 
In part two,  we move the roots $\{\alpha_{h,l}\}_{l=h+1}^{l=N}\subset A_h$
to the right through the set $B_{h-1}$. The result after part two is the order $\succ_{h-1}$. 

\begin{figure}[h]
\begin{center}
\scalebox{.8}{\includegraphics{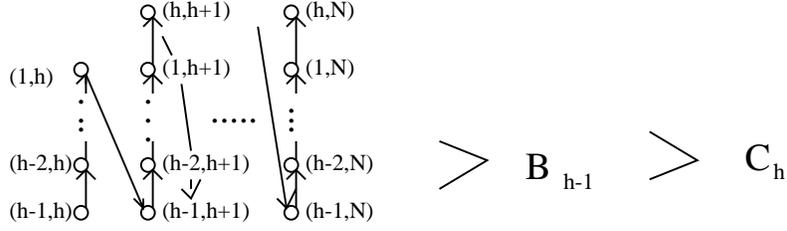}}
\end{center}
\caption{The resulting order after part one}\label{Fi:part1}
\end{figure}

Part one. The root $\alpha_{h-1,h}$ is the largest root in $B_h$. Thus, it
immediately succeeds the set of roots $A_h$ with respect to the order 
$\succ_h$. On the other hand, it is the highest element of the order $\succ_{h-1}$. We need to move
$\alpha_{h-1,h}$ through all elements of the set $A_h$. Make $h-2$  elementary 
transformations of type $A_1\oplus A_1$ moving $\alpha_{h-1,h}$ over 
$\alpha_{1,N},\ldots,\alpha_{h-2,N}$. Then perform
 $\cdots, \alpha_{h,N},\alpha_{h-1,N},\alpha_{h-1,h},\cdots\rightarrow$ 
$\cdots, \alpha_{h-1,h},\alpha_{h-1,N},\alpha_{h,N},\cdots$. Again, make $h-2$  elementary 
transformations of type $A_1\oplus A_1$ moving $\alpha_{h-1,h}$ 
over $\alpha_{1,N-1},\ldots,\alpha_{h-2,N-1}$. Then perform
$\cdots, \alpha_{h,N-1},\alpha_{h-1,N-1},\alpha_{h-1,h},\cdots \rightarrow $
$\cdots, \alpha_{h-1,h},\alpha_{h-1,N-1},\alpha_{h,N-1},\cdots$. After $N-h$ groups of steps of 
the above type the root $\alpha_{h-1,h}$ is greater than any root of the set $A_h$. 
Now, the root $\alpha_{h-2,h}$ is the largest root in $B_h-\{\alpha_{h-1,h}\}$ and it 
immediately succeeds the set $A_h$ in the normal order resulting after the above $N-h$ groups of steps. 
Repeat the procedure
for the root $\alpha_{h-2,h}$, and then for the roots $\alpha_{h-3,h},\ldots,\alpha_{1,h}$. 
This is the end of part one. The resulting order after part one is shown in Figure~\ref{Fi:part1}. 
In part one we used exactly one $A_2$ elementary transformation of  type
$\cdots, \alpha_{h,l},\alpha_{k,l},\alpha_{k,h},\cdots\rightarrow 
\cdots, \alpha_{k,h},\alpha_{k,l},\alpha_{h,l},\cdots$
for every pair of positive integers $(k,l)$, $1\leq k<h<l\leq N$.
All other elementary transformations in part one were  of the type  $A_1\oplus A_1$.

Part two. We need to move the roots $\alpha_{h,N},\ldots,\alpha_{h,h+1}$ from $A_h$ through $B_{h-1}$ 
to put them in front of $C_h$, see Figure~\ref{Fi:part1}. We first move
$\alpha_{h,N}$, then $\alpha_{h,N-1}$ and 
so on. It is easy to see that all elementary transformations we use in this part are of type 
$A_1\oplus A_1$. \qed

\subsection{Intermediate normal orders between $\succ_{h}$ and $\succ_{h-1}$.}\label{HGsol_p}
Consider again the sequence of elementary transformations $\sigma=\{\sigma_1,\ldots,\sigma_M\}$.
Denote $\succ^{p}$ the normal order we obtain from $\succ_h$ after applying $\sigma_1,\ldots,\sigma_p$.
Set $\succ^{0}=\succ_{h}$. We have $\succ^{M}=\succ_{h-1}$.
From the construction of the sequence $\sigma$ it follows that  
$A_2$ elementary transformations in $\sigma$ are labeled by  pairs 
$(k,l)$,  $1\leq k<h<l\leq N$. Namely, a pair $(k,l)$ labels the transformation 
$$\cdots,\alpha_{h,l},\alpha_{k,l},\alpha_{k,h},\cdots\rightarrow
\cdots,\alpha_{k,h},\alpha_{k,l},\alpha_{h,l},\cdots.$$ For $p=0,\ldots,M$, let
$X_p$ be the set of pairs $(k,l)$ which label  $A_2$ transformations in 
$\{\sigma_1,\ldots,\sigma_p\}$. 

For every pair of integers $(k,l)$, $k<l$, set $ F_{k,l}(\succ^p)=F_{k,l}(h-1)$ for $(k,l)\in X_p$, and
$F_{k,l}(\succ^p)=F_{k,l}(h)$ for $(k,l)\not\in X_p$.
The set $\{-F_{k,l}(\succ^p)\}_{k<l}$ is a basis of $\gnm$. Order it according to $\succ^p$. 
The corresponding PBW-basis of $U(\gnm)$ is 
$$F_I(\succ^p)= (-1)^{\sum i_{l,k}}\prod_{k<l}\frac{F_{k,l}(\succ^p)^{i_{l,k}}}{i_{l,k}!}.$$
Let $\{F_{I_1}(\succ^p),\ldots,F_{I_{n+1}}(\succ^p)\}$ be elements of the PBW-basis of
$U(\gnm)$ corresponding to $\succ^p$. Set $I=(I_1,\cdots,I_{n+1})$ and 
$F_I(\succ^p)v=F_{I_1}(\succ^p)v_1\otimes\cdots\otimes F_{I_{n+1}}(\succ^p)v_{n+1}$.
The vector $F_I(\succ^p)v$ belongs to $V'[\nu']$ if $I\in P(\bm,n+1)$
 
Introduce the rational function $\phi(I;\succ^{p})$ as
$$\phi(I;\succ^{p})=\prod_{j=1}^{n+1}\left(\prod_{(k',l')\in X_p}\prod_{q=1}^{i_{l',k'}}
\phi_{h-1}^{(j,k',l',q)}\times
\prod_{(k,l)\not\in X_p}\prod_{q=1}^{i_{l,k}}\phi_h^{(j,k,l,q)}\right).$$

\subsection{ Geometric interpretation of the rational functions $\phi(I,h)$ and $\phi(I;\succ^{p})$.}
\label{HGgeom}
\begin{figure}[h]
\begin{center}
\scalebox{.7}{\includegraphics{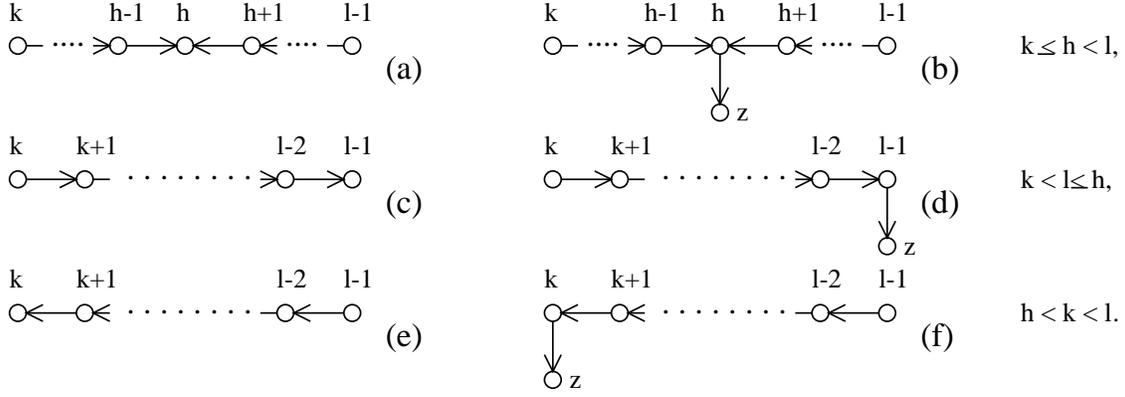}}
\end{center}
\caption{Strings and grounded strings corresponding to $\succ_h$}\label{Fi:string_j}
\end{figure}

For any pair of integers $(l,k)$, $1\leq k<l\leq N$, let the string of type $(l,k)$ 
corresponding to $\succ_h$ be the oriented graph with $l-k$ vertices labeled 
by positive integers $k,\ldots,l-1$ given in Figure~\ref{Fi:string_j}, cases (a), (c), (e).
For any $z\in\cit$, 
let the string of type $(l,k)$ corresponding to  $\succ_h$ grounded at $z$ be the oriented tree
with labeled vertices given in Figure~\ref{Fi:string_j}, cases (b), (d), (f).

\begin{figure}[h]
\begin{center}
\scalebox{.7}{\includegraphics{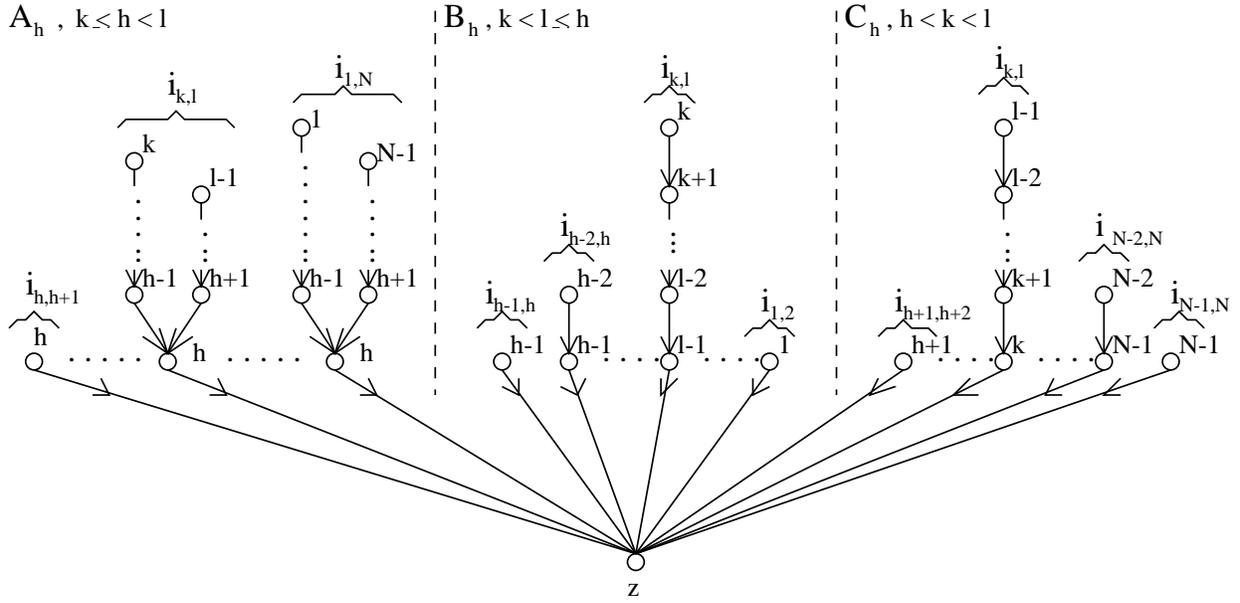}}
\end{center}
\caption{The tree $T(I_0,h,z)$}\label{Fi:tree_j}
\end{figure}

Let $I_0=\{i_{l,k}\}_{k<l}$ be a set of non-negative integers.
Define the tree $T(I_0,h,z)$ as the union of strings  corresponding to $\succ_h$ 
grounded at $z$, such that the string of 
type $(l,k)$ grounded at $z$ enters exactly $i_{l,k}$ times and all vertices $z$ of all strings 
are then identified, see Figure~\ref{Fi:tree_j}.
Define the tree $T(I_0,\succ^p,z)$ as follows.
For every $(k,l)\in X_p$ take $i_{l,k}$  strings of 
type $(l,k)$ corresponding to  $\succ_{h-1}$ and grounded at $z_1$.  
For every $(k,l)\not\in X_p$ take $i_{l,k}$ strings of 
type $(l,k)$  corresponding to  $\succ_{h}$ and grounded at $z_1$. 
Let $T(I_0,\succ^p,z)$ be the union of all selected strings with all vertices $z_1$  then identified.

For any $I=(I_1,\ldots I_{n+1})\in P(\bm,n+1)$, let the forest $T(I,h)=\sqcup_{j=1}^{n+1}T(I_j,h,z_j)$ 
be the disjoint union of  trees, and let the forest $T(I,\succ^p)=\sqcup_{j=1}^{n+1}T(I_j,h,\succ^p)$ 
be the disjoint union of  trees. 

For any $p=1,\ldots, N-1$, the number of vertices in the forest 
$T(I,h)$ labeled by $p$ equals $m_p$. For every $p$ fix a bijection $\beta^{h}_p(I):$ $\{$the vertices of
$T(I,h)$ labeled by $p\}$ $\rightarrow S_p$. To every oriented edge $(w',w'')$
from  a vertex $w'$ labeled by an integer $p'$ to  a vertex $w''$ labeled by an integer $p''$ assign  
the function $f(w',w'')=1/(t_{p'}^{(\beta^{h}_{p'}(I)(w'))}-t_{p''}^{(\beta^{h}_{p'}(I)(w''))})$. 
To every oriented edge $(w',w'')$
from  vertex $w'$ labeled by an integer $p'$ to  a vertex $w''$ labeled by a complex number $z$ assign  
the function $f(w',w'')=1/(t_{p'}^{(\beta^{h}_{p'}(I)(w'))}-z)$. 
Define  rational functions corresponding to a forest (respectively,
a tree, a grounded string, or a string) as the product of the functions $f(w',w'')$ over all edges $(w',w'')$
of the  forest (respectively, the tree, the grounded string, or the string). 
It is clear that the rational function corresponding to the forest $T(I,h)$ is exactly the function
$\phi(I,h)$ defined in Section~\ref{spec_ord}.

Define a rational function corresponding to the forest $T(I,\succ^p)$ analogously. 
The rational function corresponding to the forest $T(I,\succ^p)$ equals $\phi(I,\succ^p)$.

\subsection{Consistency of the definitions.} \label{consistency}
Let $h'\in\{h-1,h\}$. Since  
$\succ^{0}=\succ_{h}$, $\succ^{M}=\succ_{h-1}$, we have given two different 
definitions of the basis of $\gnm$ corresponding to $\succ_{h'}$, the rational function $\phi(I,h')$, 
and the forest $T(I,h')$, see Sections~\ref{HGsol_h}, \ref{HGsol_p}, \ref{HGgeom}. 
We will show that these two definitions are consistent.

  In the case $h'=h$, we have $X_0=\varnothing$. Trivially, 
$T(I,h)=T(I,\succ^0)$, $\phi(I,h)=\phi(I,\succ^0)$, and $F_{k,l}(h)=F_{k,l}(\succ^0)$ for $k<l$.

  In the case $h'=h$, we have $X_M=\{(k,l)\,|\,1\leq k<h<l\leq N\}$.  We start with the  basis of $\gnm$.
Recall that $F_{k,l}(h')=(-1)^{a_{k,l}(h')}e_{l,k}$, for any $h'$. 
The following five cases exhaust all possibilities.\\
$\bullet$ If $k<h<l$, then  $F_{k,l}(\succ^M)=F_{k,l}(h-1)$ by definition.\\
$\bullet$ If $k<l< h$, then $F_{k,l}(\succ^M)=F_{k,l}(h)$, and $a_{k,l}(h)=a_{k,l}(h-1)=0$. Thus,
$F_{k,l}(h)=F_{k,l}(h-1)$ and $F_{k,l}(\succ^M)=F_{k,l}(h-1)$.\\
$\bullet$ If $k<l= h$, then $F_{k,l}(\succ^M)=F_{k,l}(h)$, $a_{k,l}(h)=0$, and $a_{k,l}(h-1)=l-(h-1)-1$. Thus,
$F_{k,l}(h)=F_{k,l}(h-1)$ and $F_{k,l}(\succ^M)=F_{k,l}(h-1)$.\\
$\bullet$ If $h<k<l$, then $F_{k,l}(\succ^M)=F_{k,l}(h)$, $a_{k,l}(h)=a_{k,l}(h-1)=l-k-1$. Thus,
$F_{k,l}(h)=F_{k,l}(h-1)$ and $F_{k,l}(\succ^M)=F_{k,l}(h-1)$.\\
$\bullet$ If $h=k<l$, then $F_{k,l}(\succ^M)=F_{k,l}(h)$, $a_{k,l}(h)=l-h-1$, and $a_{k,l}(h-1)=l-k-1$. Thus,
$F_{k,l}(h)=F_{k,l}(h-1)$ and $F_{k,l}(\succ^M)=F_{k,l}(h-1)$.\\

 Similarly, for every set of non-negative integers $I_0=\{i_{l,k}\}_{k<l}$, 
we compare the tree $T(I_0,z,\succ^M)$ with the tree $T(I_0,z,h-1)$.
In order to construct the tree $T(I_0,z,h-1)$, for every pair $(l,k)$, $k<l$, we took $i_{l,k}$ strings of type
$(l,k)$ grounded at $z_1$ which correspond to the order $\succ_{h-1}$ and then we identified all vertices $z_1$.
 In the construction of the tree $T(I_0,z,\succ^M)$,  for every pair $(l,k)$, $k<h<l$,  
we took $i_{l,k}$ strings of type $(l,k)$ grounded at $z_1$ which correspond to the order $\succ_{h-1}$, and
for every pair $(l,k)$, $h\leq k<l$ or $k<l\leq h$,  
we took $i_{l,k}$ strings of type $(l,k)$ grounded at $z_1$ which correspond to the order $\succ_{h}$.
Then we identified all vertices $z_1$. 
 For every pair $(l,k)$, $h\leq k<l$ or $k<l\leq h$,  the string of type $(l,k)$ 
grounded at $z_1$ which corresponds to the order $\succ_{h-1}$ coincides with the string of type $(l,k)$ 
grounded at $z_1$ which corresponds to the order $\succ_{h}$, see Figure~\ref{Fi:string_j}. Therefore, 
the tree $T(I_0,z,\succ^M)$ coincides with the tree $T(I_0,z,h-1)$.
 Conclude that the forest $T(I,h-1)$ coincides with the forest $T(I,\succ^M)$, and 
we have $\phi(I,h-1)=\phi(I,\succ^M)$.

\subsection{Proof of Theorem~\ref{mat_trig_j}.} The theorem states that, for any $h$, 
we have  $\phi(z,t;h)=\phi(z,t)$.
The equality $\phi(z,t;N-1)=\phi(z,t)$  is valid by definition. 
It suffices to show that $\phi(z,t;h)=\phi(z,t;h-1)$
for every $h$.
 
First consider the case $n=0$. Recall that $n+1$ is the number of tensor factors in $V'$. Thus,
$V'=V_1$ is a Verma module. The discussion in Section~\ref{consistency} implies that
$$\phi(z,t;h)=\sum_{I\in P(\bm,1)}\phi(I,\succ^0)F_I(\succ^0)v_1,\qquad
\sum_{I\in P(\bm,1)}\phi(I,\succ^M)F_I(\succ^M)v_1=\phi(z,t;h-1).$$
In order to finish the case $n=0$, it suffice to show that
\begin{equation}\label{statementB}
\sum_{I\in P(\bm,1)}\phi(I,\succ^p)F_I(\succ^p)v_1
=\sum_{I\in P(\bm,1)}\phi(I,\succ^{p+1})F_I(\succ^{p+1})v_1, \mbox{ for } p=0,\ldots,M-1.
\end{equation}

  Fix $p\in\{0,\ldots,M-1\}$.
\paragraph{} Assume that $\sigma_{p+1}$ is of type $A_1\oplus A_1$.
 We have $X_p=X_{p+1}$ and $\phi(I,\succ^p)=\phi(I,\succ^{p+1})$ for every index $I\in P(\bm,1)$.
The case $A_1\oplus A_1$ of Lemma~\ref{elem_trans} implies $F_I(\succ^p)=F_I(\succ^{p+1})$ 
for every index $I\in P(\bm,1)$. This implies (\ref{statementB}).
\paragraph{} Assume that $\sigma_{p+1}$ is of type $A_2$ and is labeled by a pair
of integers $(k,l)$, $k<h<l$. The set $X_{p+1}$ equals the set $X_p\cup\{(k,l)\}$.
Let us check, that we can apply the $A_2$ case of Lemma~\ref{elem_trans} to 
the triple $F_{h,l}(\succ^p),F_{k,l}(\succ^p),F_{k,h}(\succ^p)\in sl_N$.
Set $f_\alpha=F_{h,l}(h)$, $f_\beta=F_{k,h}(h)$ . Then $f_\alpha=F_{h,l}(\succ^p)=F_{h,l}(\succ^{p+1})$,
$f_\beta=F_{k,h}(\succ^p)=F_{k,h}(\succ^{p+1})$, 
\begin{align} [f_\alpha,f_\beta]&=(-1)^{l-h-1}[e_{l,h},e_{h,k}]=
(-1)^{l-h-1}e_{l,k}=F_{k,l}(h)=F_{k,l}(\succ^p),\notag\\
[f_\beta,f_\alpha]&=(-1)^{l-(h-1)-1}e_{l,k}= F_{k,l}(h-1)=F_{k,l}(\succ^{p+1}). \notag
\end{align}
Here we used the explicit form of the numbers $\{a_{k,l}(h)\}$, namely,
$$ a_{k,l}(h)=0 \mbox{ if } k<l\leq h,\quad a_{k,l}(h)=l-k-1 \mbox{ if } h<k<l,\quad
a_{k,l}(h)=l-h-1 \mbox{ if } k\leq h<l.$$
By Lemma~\ref{elem_trans}, we have
\begin{equation}\notag
 \frac{F_{h,l}(\succ^p)^a}{a!}\frac{F_{k,l}(\succ^p)^c}{c!}\frac{F_{k,h}(\succ^p)^b}{b!}=
\sum_r \bin{c+r}{r}(-1)^{c-r}\frac{F_{k,h}(\succ^{p+1})^{b-r}}{(b-r)!}
\frac{F_{k,l}(\succ^{p+1})^{c+r}}{(c+r)!}\frac{F_{h,l}(\succ^{p+1})^{a-r}}{(a-r)!}.
\end{equation}
Use this equality to transform $\sum_{I\in P(\bm,1)}\phi(I,\succ^p)F_I(\succ^p)v_1$. 
Set $a=i_{h,l}$, $b=i_{k,h}$, $c=i_{k,l}$,
$I'=\{i_{p,q}\}_{(p,q)\not\in\{(h,l),(k,h),(k,l)\}}$, $|I'|=\sum i_{p,q}$. We have
$I=\{I',a,b,c\}$ and
\begin{align}
\sum_{I\in P(\bm,1)}&\phi(I,\succ^p)F_I(\succ^p)v_1=\sum_{I',d_1,d_2}(-1)^{|I'|}
\sum_{ \begin{smallmatrix} a+c=d_1,\\b+c=d_2\end{smallmatrix}}
(-1)^{a+b+c}\phi(\{I',a,b,c\};\sw)\times\notag\\
&\times\cdots\frac{F_{h,l}(\succ^p)^a}{a!}\frac{F_{k,l}(\succ^p)^c}{c!}
\frac{F_{k,h}(\succ^p)^b}{b!}\cdots v_1\notag\\
&=\sum_{I',d_1,d_2}(-1)^{|I'|}\sum_{ \begin{smallmatrix} a+c=d_1,\\b+c=d_2\end{smallmatrix}}
\sum_r \bin{c+r}{r}(-1)^{c}(-1)^{a+b+c-r}\phi(\{I',a,b,c\};\sw) \times\notag\\
&\times\cdots\frac{F_{k,h}(\succ^{p+1})^{b-r}}{(b-r)!}
\frac{F_{k,l}(\succ^{p+1})^{c+r}}{(c+r)!}\frac{F_{h,l}(\succ^{p+1})^{a-r}}{(a-r)!}\cdots v_1.\notag
\end{align}
\begin{align}
\sum_{I\in P(\bm,1)}&\phi(I,\succ^p)F_I(\succ^p)v_1
=\sum_{I',d_1,d_2}(-1)^{|I'|}\sum_{ \begin{smallmatrix} a'+c'=d_1,\\b'+c'=d_2\end{smallmatrix}}
(-1)^{a'+b'+c'}\psi(\{I',a',b',c'\};\succ^{p})\times\notag\\
\label{eq40}
&\times\cdots\frac{F_{k,h}(\succ^{p+1})^{b'}}{(b')!}
\frac{F_{k,l}(\succ^{p+1})^{c'}}{(c')!}\frac{F_{h,l}(\succ^{p+1})^{a'}}{(a')!}\cdots v_1,
\end{align}
where $\psi(\{I',a',b',c'\};\succ^{p})= \sum_r \bin{c'}{r}(-1)^{c'-r}\phi(\{I',a'+r,b'+r,c'-r\};\sw)$.

The following claim implies (\ref{statementB}).\\\hfill\\
{\bf Claim.} For every $I=\{I',a',b',c'\}\in P(\bm,1)$ we have 
\begin{equation}\label{claim}
\psi(\{I',a',b',c'\};\succ^{p})=\phi(\{I',a',b',c'\};\succ^{p+1}).
\end{equation}
\paragraph{} The right hand side of (\ref{claim}) equals the product of the functions $f(w',w'')$ over
all oriented edges $(w',w'')$ of the tree $T(I,z_1,\succ^{p+1})$. The left hand side of (\ref{claim})
depends on rational functions $\{\phi(\{I',a'+r,b'+r,c'-r\};\sw)\}_r$ which correspond to
trees $\{T(\{I',a'+r,b'+r,c'-r\},z_1,\succ^{p+1})\}_r$. For every $r$,
the tree $T(I,z_1,\succ^{p+1})$
and the tree $T(\{I',a'+r,b'+r,c'-r\},z_1,\succ^{p+1})$  differ only at the grounded
strings of types $(k,l),(h,l),(k,h)$. Namely, the tree $T(I,z_1,\succ^{p+1})$ has
$a'$ strings of type $(h,l)$ grounded at $z_1$ corresponding to $\succ_h$,
$b'$ strings of type $(k,h)$ grounded at $z_1$ corresponding to $\succ_h$ and
$c'$ strings of type $(k,l)$ grounded at $z_1$ corresponding to $\succ_{h-1}$.
The tree $T((I',a'+r,b'+r,c'-r),z_1,\succ^{p+1})$ has
$a'+r$ strings of type $(h,l)$ grounded at $z_1$ corresponding to $\succ_h$,
$b'+r$ strings of type $(k,h)$ grounded at $z_1$ corresponding to $\succ_h$ and
$c'-r$ strings of type $(k,l)$ grounded at $z_1$ corresponding to $\succ_h$.

\begin{figure}[h]
\begin{center}
\scalebox{.7}{\includegraphics{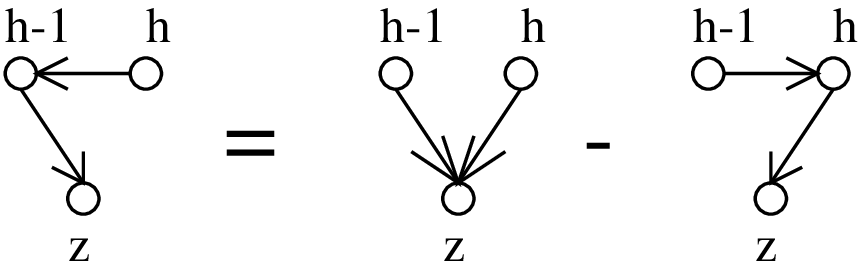}}
\end{center}
\caption{}\label{Fi:trngl}
\end{figure}

\begin{figure}[h]
\begin{center}
\scalebox{.7}{\includegraphics{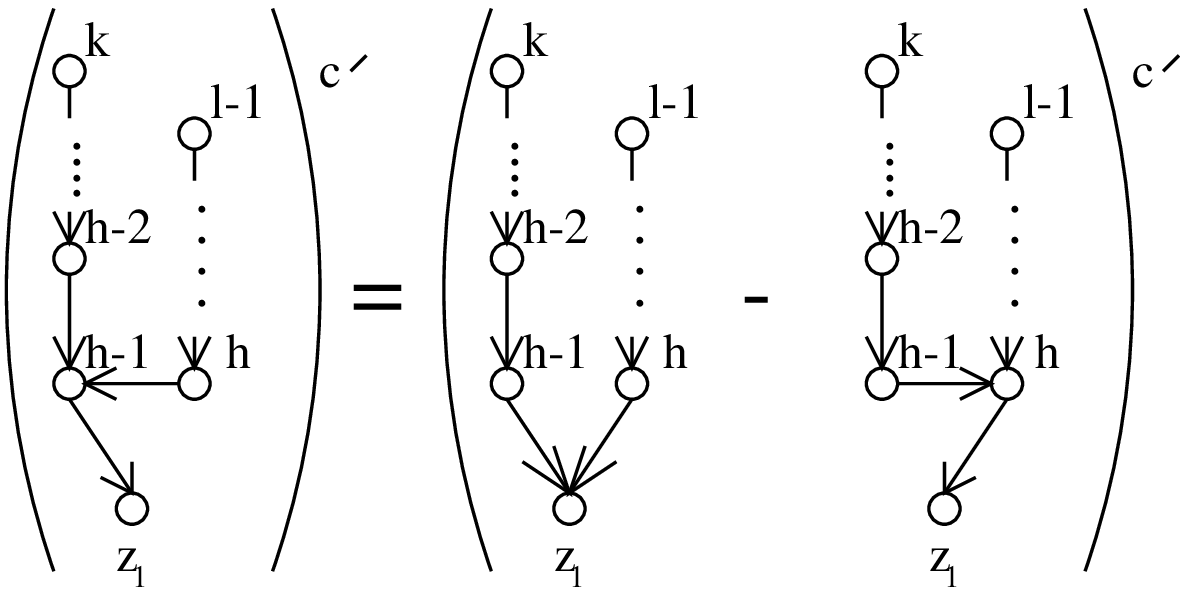}}
\end{center}
\caption{}\label{Fi:j-1_to_j}
\end{figure}

\begin{figure}[h]
\begin{center}
\scalebox{.7}{\includegraphics{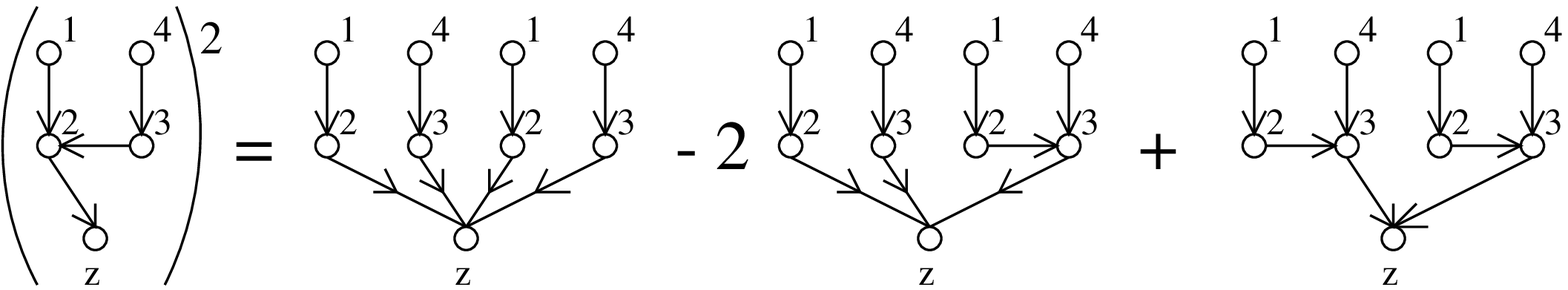}}
\end{center}
\caption{Example}\label{Fi:ex}
\end{figure}

The claim now follows from identities
\begin{align}\label{trng}
\frac{1}{(t_h-t_{h-1})(t_{h-1}-z)}&=\frac{1}{(t_h-z)(t_{h-1}-z)}-\frac{1}{(t_{h-1}-t_h)(t_{h}-z)},\\
\prod_{q=1}^{c'}\phi^{(1,k,l,q)}_{h-1}&=\prod_{q=1}^{c'}(\phi^{(1,k,h,q)}_{h}\phi^{(1,h,l,q)}_{h}-
\phi^{(1,k,l,q)}_h)\label{tr_pow}
\end{align}
Identity (\ref{trng}) has a simple tree interpretation shown in Figure~\ref{Fi:trngl}.
\\
Identity (\ref{tr_pow}) has a tree interpretation shown in Figure~\ref{Fi:j-1_to_j}.
In Figure~\ref{Fi:j-1_to_j} multiplication of trees means  disjoint union of trees,  addition is formal,
and multiplication distributes through addition, see the example in Figure~\ref{Fi:ex}. 
The claim is proved. Thus, the statement of Theorem~\ref{mat_trig_j} is proved for $n=0$. \qed

The statement of Theorem~\ref{mat_trig_j} for arbitrary $n$ is proved as the case $n=0$ since 
we execute the changes of bases and the reorganization of rational functions 
separately for each tensor factor.
There is no interference between the rational functions corresponding to different tensor factors 
since the grounding points of the strings are distinct. \qed

\section{Proof of Theorem~\ref{main_sln}.}\label{mainproof}
Let $V=M_1\otimes\cdots\otimes M_n$, where $M_j$ is the $sl_N$ Verma module 
with highest weight $\Lambda_j$ and highest weight vector $v_j$. Set
$v=v_1\otimes\cdots\otimes v_n$.  
Fix a weight subspace $V[\nu]\subset V$, where
$\nu=\sum_{j=1}^n\Lambda_j-\sum_{k=1}^{N-1}m_k\alpha_k$, and $\sum_{k=1}^{N-1}m_k\alpha_k\in Q^+$. 
Recall the notation $m=\sum_{k=1}^{N-1}m_k$, $\bm=\sum_{k=1}^{N-1}m_k\alpha_k$. Set 
$\Lambda_{n+1}=\lambda-\rho-\frac12\nu$.
Let $M_{n+1}$ be the $sl_N$ Verma module 
with highest weight $\Lambda_{n+1}$ and highest weight vector $v_{n+1}$. 
Consider the auxiliary space $V'=V\otimes M_{n+1}$ 
and its weight space $V'[\nu']$, where $\nu'=\nu+\Lambda_{n+1}$.
 
\subsection{Vector spaces of rational functions}. 
For every $I\in P(\bm,n+1)$ fix a set of bijections $\{\beta_p(I)\}$ and define 
a vector space of rational functions 
$\Omega'(\{t_k^{(d)}\},\{z_j\})=\cit\{\phi(I)\,|\, I\in P(\bm,n+1)\}$. These functions depend on the
complex variables $t_k^{(d)},z_j$, $k=1,\ldots,N-1$, $d\in\{1,\ldots,m_k\}$, $j=1,\ldots,n+1$. 
Denote $\Omega'$ the restriction $\Omega'(\{t_k^{(d)}\},z_1,\ldots,z_n,0)$. We will use the same notation
$\phi(I)$ for a function in $\Omega'$.

For every $I\in P(\bm,n)$ fix a set of bijections $\{\beta_p(I)\}$ and define 
a vector space of rational functions $\Omega=\cit\{\phi(I)\,|\, I\in P(\bm,n)\}$.
These functions depend on the
complex variables $t_k^{(d)},z_j$, $k=1,\ldots,N-1$, $d\in\{1,\ldots,m_k\}$, $j=1,\ldots,n$.
Identify each index $I\in P(\bm,n)$ with the index $I'=(I,\{i_{l,k}^{(n+1)}=0\}_{k<l})\in P(\bm,n+1)$.
Since $\phi(I')=\phi(I)$ for an appropriate choice of bijections $\{\beta_p(I')\}$, $\{\beta_p(I)\}$ 
we will consider $\Omega$ as a subset of $\Omega'$.

Consider the set of rational functions with poles in $\D$, see Section~\ref{HGsol}. Note that it contains 
$\Omega'$. Define an equivalence relation $\backsim$ on the set. Namely, let $\phi_1,\phi_2$ be 
two elements of the set.
We write $\phi_1\backsim\phi_2$ if $\int_{\gamma(z)}\Phi^{1/\kappa}\phi_1\,dt=
\int_{\gamma(z)}\Phi^{1/\kappa}\phi_2\,dt$ for any horizontal family of integration cycles $\gamma(z)$, 
which satisfies Matsuo's assumption. 

\paragraph{\bf Examples.} (a) For a fixed $I\in P(\bm,n+1)$, the functions $\phi(I)$ corresponding
to different choices of bijections $\{\beta_p(I)\}$ are equivalent under $\backsim$. \\
(b) If there exists a rational form $\eta$ with poles in $\D$ such that $\Phi^{\frac1{\kappa}}
(\phi_1-\phi_2)dt=d_t(\Phi^{\frac1{\kappa}}\eta)$, 
where $d_t$ denotes the exterior differentiation with respect to the variables $t_k^{(d)}$,
then $\phi_1\backsim\phi_2$ because $\gamma(z)$ is a family of closed cycles.

The set $\{(F_Iv')^*\}_{I\in P(\bm,n+1)}$ is a basis of $(V'[\nu'])^*$, see the definition of $\{F_Iv'\}$ 
in Section~\ref{HGsol}. Define a linear map $D':(V'[\nu'])^*\rightarrow \Omega'$
by setting $D((F_Iv')^*)=\phi(I)$, for all $I\in P(\bm,n+1)$. 
Analogously, define a linear map $D:(V[\nu])^*\rightarrow \Omega$
by setting $D((F_Iv)^*)=\phi(I)$, for all $I\in P(\bm,n)$.

For indices $I\in P(\bm,n)$ and $I'=(I,\{i_{l,k}^{(n+1)}=0\}_{k<l})\in P(\bm,n+1)$ we have
\begin{equation}\label{D'=D}
D'(F_{I'}v')=D'(F_Iv\otimes v_{n+1})=\phi(I')=\phi(I)=D(F_Iv).
\end{equation}

Reformulate Theorem~\ref{mat_trig_j} as follows.
\begin{theorem} \label{mat_trig_j1}For every $h=1,\ldots,N-1$, we have
  $\sum_{I\in P(\bm,n)} D((F_I(h)v)^*) F_I(h)v=\phi(t,z)$. Moreover, for all $I$, 
$D((F_I(h)v)^*)=\phi(I,h)$.
\end{theorem}

\subsection{The action of $sl_N$ on $(V')^*$.} 
For every set of non-negative integers $I_0=\{i_{l,k}\}_{k<l}$, and every $l',k'$, $k'<l'$, denote
$I_0\pm 1_{l',k'}=\{i'_{l,k}\}_{k<l}$, where $i'_{l',k'}=i_{l',k'}\pm 1$ and $i'_{l,k}=i_{l,k}$ if
$(l,k)\ne (l',k')$. Let $I=(I_1,\ldots,I_{n+1})$, where $I_j=\{i_{l,k}^{(j)}\}_{k<l}$ 
is a set of non-negative integers.  Let $I\pm 1_{l',k'}^{(j)}$ denote addition in the component $I_j$.
\begin{lemma}[cf. Lemma 4.2 in \cite{Mat}] \label{E-act}
For every $I$, and every $h=1,\ldots,N-1$, we have 
\begin{align}
E_{\alpha_h}(F_Iv')^* &=\sum_{j=1}^{n+1}\left( \sum_{p=h+2}^{N-1}i_{p,h+1}^{(j)}
(F_{I+1_{p,h}^{(j)}-1_{p,h+1}^{(j)}}v')^*
-\sum_{p=1}^{h-1}i_{h,p}^{(j)}(F_{I+1_{h+1,p}^{(j)}-1_{h,p}^{(j)}}v')^* +\right. \notag\\
& \left. + ((\Lambda_j,\alpha_h) +\sum_{p=1}^{h-1}i_{h,p}^{(j)}-
\sum_{p=1}^h i_{h+1,p}^{(j)})(F_{I+1_{h+1,h}^{(j)}}v')^*\right).\notag
\end{align}
\end{lemma}
The proof is a straightforward computation. We use the standard PBW-basis of $U(\gnm)$.\qed
\begin{lemma}\label{F-act}
For every $I$, and every $k,h,l$, $1\leq k\leq h <l\leq N$, we have
\begin{equation}
 F_{k,l}(h)(F_I(h)v')^*=\sum_{j=1}^{n+1}i_{l,k}^{(j)}(F_{I-1_{l,k}^{(j)}}(h)v')^*.
\end{equation}
\end{lemma}
The proof is a straightforward computation. We use the PBW-basis of $U(\gnm)$ corresponding to $\succ_h$.\qed

  Next, we will compute the image of $E_{\alpha_h}(F_Iv')^*\in V'[\nu']$ under $D'$. According to 
Lemma~\ref{E-act}  $D'(E_{\alpha_h}(F_Iv')^*)$ is a sum of rational functions with complex coefficients,
$\sum_{J\in P(\bm,n+1)} c_J\phi(J)$. For each index $J$, we can choose a set of bijections 
$\{\beta_p(J)\}$ such that the function $\phi(J)$ has the form $\frac1{t_h-*}\phi(I)$, where 
$t_h$ is a new variable of type $h$, and $*$ is any of the $t$ variables fixed by the set of bijections
$\{\beta_p(I)\}$, or a complex number in the set $\{z_1,\ldots,z_n,0\}$.

\begin{lemma}[\cite{Mat}, Lemma~3.4]\label{mat-l.3.4} Let $s=(j,k,l,q)\in S(I)$. Then, we have
\begin{align}
&(a)\,\,\frac1{t_h-z_j}\phi^{(j,k,l,q)}\backsim
 \phi^{(j,h,h+1,i^{(j)}_{h+1,h}+1)}\phi^{(j,k,l,q)},\notag\\
&(b)\,\,\left(\sum_{p=k}^{l-1}\frac{(\alpha_h,\alpha_p)}{t_h-t_p^{(\beta_p(I)(s))}}\right)
\phi^{(j,k,l,q)}
\backsim\left\{ \begin{matrix}
-\phi^{(j,h,l,i^{(j)}_{l,h}+1)} & \mathrm{if} & h=k-1,\\
\phi^{(j,k,l,q)}\phi^{(j,h,h+1,i^{(j)}_{h+1,h}+1)} &\mathrm{if} & h=l-1,\\
\phi^{(j,k,h+1,i^{(j)}_{h+1,l}+1)}-\phi^{(j,k,h,q)}\phi^{(j,h,h+1,i^{(j)}_{h+1,h}+1)} & 
\mathrm{if} & h=l,\\
0 & &\mathrm{otherwise}.
\end{matrix}\right.\notag
\end{align}
\end{lemma}

\begin{corollary} 
\begin{align}
(a)\quad\frac{d_{t_h}\Phi}{\Phi}\phi(I)&\backsim\sum_{j=1}^{n+1}\left( -\sum_{p=h+2}^{N-1}i_{p,h+1}^{(j)}
\phi(I+1_{p,h}^{(j)}-1_{p,h+1}^{(j)})
+\sum_{p=1}^{h-1}i_{h,p}^{(j)}\phi(I+1_{h+1,p}^{(j)}-1_{h,p}^{(j)}) +\right. \notag\\
& \left. + (-(\Lambda_j,\alpha_h) -\sum_{p=1}^{h-1}i_{h,p}^{(j)}+
\sum_{p=1}^h i_{h+1,p}^{(j)})\phi(I+1_{h+1,h}^{(j)})\right).\notag
\end{align}

$(b)$ $D'(E_{\alpha_h}(F_Iv')^*) \backsim \frac{d_{t_h}\Phi}{\Phi}\phi(I)$.
\end{corollary}

{\it Proof.} (a) Our choice of $\Lambda_{n+1}$ is consistent with the definition of the function $\Phi$,
that is $(t_k^{(d)})^{-(\alpha_k,\Lambda_{n+1})}$ is the contribution of the point $z_{n+1}=0$ 
in $\Phi$.
Observe that
$$\frac{d_{t_h}\Phi}{\Phi}=\sum_{j=1}^{n+1}\left(
-(\Lambda_j,\alpha_h)\frac{d(t_h-z_j)}{(t_h-z_j)}+\sum_{k<l}\sum_{q=1}^{i_{l,k}^{(j)}}
\sum_{h=k}^{l-1}(\alpha_h,\alpha_p)\frac{d(t_h-t_p^{(\beta_p(I)(s))})}{(t_h-t_p^{(\beta_p(I)(s))})}\right) ,$$
and apply Lemma~\ref{mat-l.3.4}. 

(b) Part (a) and Lemma~\ref{E-act} imply statement (b). \hfill $\qed$

\begin{corollary}\label{switch-cor}
 For any $x\in U(\gnp)$ and any $u'=u\otimes u_{n+1} \in V'=V\otimes V_\lambda$, such that $x u\in V'[\nu']$,  
we have $D'(u\otimes x u_{n+1})\backsim  D'(A(x) u\otimes u_{n+1})$, where $A$ is the antipode map. 
\end{corollary}
{\it Proof.} It is sufficient to show that for any $h=1,\ldots,N-1$, we have 
$D'(E_{\alpha_h}(F_Iv')^*)\backsim 0$. Equivalently, we may show $\frac{d_{t_h}\Phi}{\Phi}\phi(I)\backsim 0$.
Introduce the logarithmic $(m-1)$-form
\begin{equation*}
\eta_I=\prod_{s\in S(I)}\left(\prod_{p=k}^{l-2}
\frac{d(t_p^{(\beta_p(I)(s))}-t_{p+1}^{(\beta_{p+1}(I)(s))})}{t_p^{(\beta_p(I)(s))}-
t_{p+1}^{(\beta_{p+1}(I)(s))}}\right)
\frac{d(t_{l-1}^{(\beta_{l-1}(I)(s))}-z_j)}{t_{l-1}^{(\beta_{l-1}(I)(s))}-z_j}.
\end{equation*}
We have $d_t\eta_I=0$, and 
 $$\kappa d_t(\Phi^{\frac1{\kappa}}\eta_I)=\Phi^{\frac1{\kappa}}\frac{d_t\Phi}{\Phi}\wedge\eta_I=
\Phi^{\frac1{\kappa}}\frac{d_{t_h}\Phi}{\Phi}\wedge\eta_I.$$
Let $i$ be the inclusion map $i:\cit^m\rightarrow \{(z,0)\}\times\cit^m\subset\cit^{n+1}\times
\cit^m$, and as usual, $dt$ denotes the standard volume form on $\cit^m_t$.
We have $i^*(\frac{d_t\Phi}{\Phi}\wedge\eta_I)=\pm\frac{d_{t_h}\Phi}{\Phi}\phi(I)dt$.
Thus, $\frac{d_{t_h}\Phi}{\Phi}\phi(I)\backsim 0$.  $\qed$

\subsection{Proof of Theorem~\ref{main_sln}.}
Let $v(z,\lambda)$ be the hypergeometric solution of the trigonometric KZ equation indicated in 
Corollary~\ref{mat_trig}. Fix $h\in\{1,\ldots,N-1\}$. 
By definition,
$$v(z,\lambda+\kappa\omega_h^\vee)= \int_{\gamma(z)}
\Phi(z,t;\lambda+\kappa\omega_h^\vee)^{1/k}\phi(z,t)\,dt.$$
The function $\Phi^{\frac{1}{\kappa}}$ changes as follows.
$$\Phi^{\frac{1}{\kappa}}(z,t;\lambda+\kappa\omega_h^\vee)=\left(
\prod_{j=1}^nz_j^{(\Lambda_j,\omega_h^\vee)}\prod_{d=1}^{m_h}\frac1{t_h^{(d)}}\right)
\Phi^{\frac{1}{\kappa}}(z,t;\lambda).$$
We study  the product $[\prod_{j=1}^nz_j^{(\Lambda_j,\omega_h^\vee)}
\prod_{d=1}^{m_h}(t_h^{(d)})^{-1}]\phi(z,t)$, using the expansion\\
$\phi(z,t)=\sum_{I\in P(\bm,n)}\phi(I,h)F_I(h)v$ obtained in Theorem~\ref{mat_trig_j}.

\begin{figure}[h]
\begin{center}
\scalebox{.7}{\includegraphics{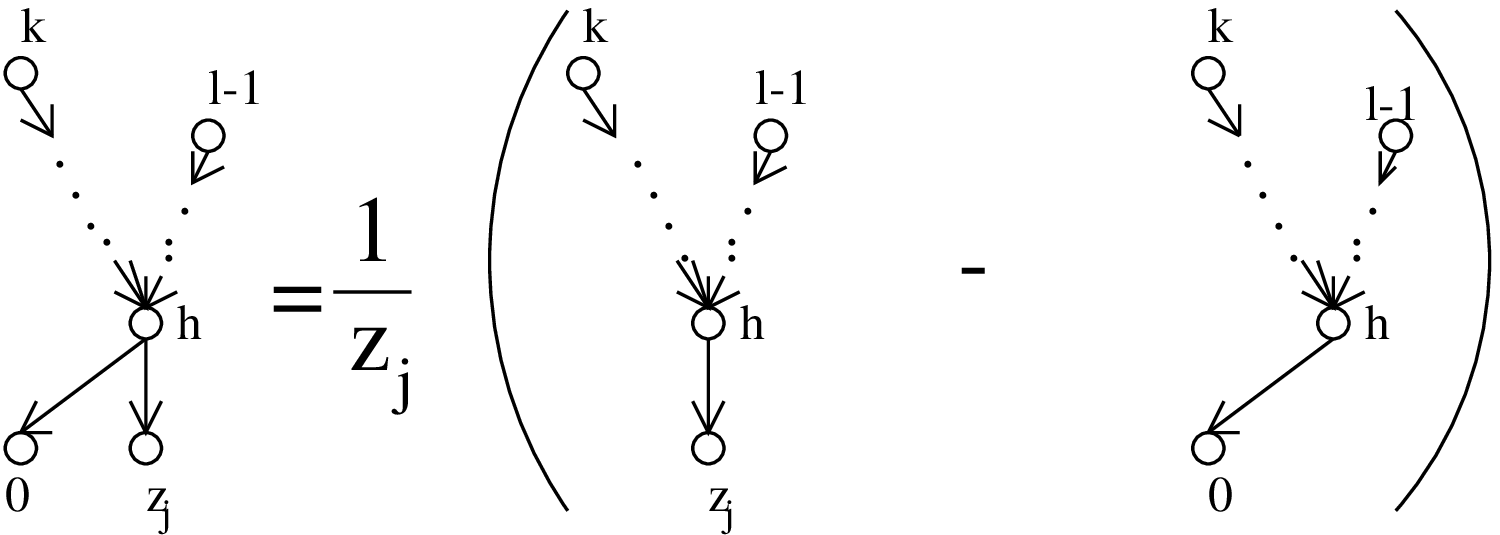}}
\end{center}
\caption{}\label{Fi:tree_0_j}
\end{figure}

The identity 
$$\frac1{(t_h^{(d)}-z_j)t_h^{(d)}}=\frac1{z_j}\left(\frac1{t_h^{(d)}-z_j}-\frac1{t_h^{(d)}}\right)$$
implies the tree interpretation shown in Figure~\ref{Fi:tree_0_j}.
Set $m_{h,j}(I)$ to be equal to the number of elements of the set $\{(j',k,l,q)\in S_h(I)\,|\,j'=j\}$. 
For all $I$, we have $\sum_{j=1}^n m_{h,j}(I)=m_h$.
Use the equality of rational functions presented in Figure~\ref{Fi:tree_0_j} and the definition
 $\phi(I,h)=\prod_{j=1}^n\prod_{(j,k,l,q)\in S(I)}\phi^{(j,k,l,q)}_h$ to obtain
\begin{equation*}
(\prod_{d=1}^{m_h}(t_h^{(d)})^{-1})\phi(I,h)=\prod_{j=1}^n z_j^{-m_{h,j}(I)}
\left( \prod_{(j,k,l,q)\in S_h(I)}(\phi^{(j,k,l,q)}_h-\phi^{(n+1,k,l,q)}_h)
\prod_{(j,k,l,q)\not\in S_h(I)}\phi^{(j,k,l,q)}_h\right).
\end{equation*}
Therefore, we have
\begin{align}[\prod_{j=1}^n & z_j^{(\Lambda_j,\omega_h^\vee)}
\prod_{d=1}^{m_h}(t_h^{(d)})^{-1}]\phi(z,t)=\label{z-factor}\\
&= \prod_{j=1}^nz_j^{(\omega_h^\vee)^{(j)}}\left( 
\prod_{(j,k,l,q)\in S_h(I)}(\phi^{(j,k,l,q)}_h-\phi^{(n+1,k,l,q)}_h)
\prod_{(j,k,l,q)\not\in S_h(I)}\phi^{(j,k,l,q)}_h\right)F_I(h)v.\notag
\end{align}
We see that the shift of the parameter $\lambda$ reduces to the following operation on  rational functions.
Take each function $\phi(I,h)$ in the formula for $\phi(z,t)$, $z\in\cit^n$, $t\in\cit^m$, 
and substitute it with the rational function
$\prod_{(j,k,l,q)\in S_h(I)}(\phi^{(j,k,l,q)}_h-\phi^{(n+1,k,l,q)}_h)
\prod_{(j,k,l,q)\not\in S_h(I)}\phi^{(j,k,l,q)}_h$, which depends on variables $(z,z_{n+1})\in\cit^{n+1}$, $t\in\cit^m$.
Then set $z_{n+1}=0$. The result is a function that depends on the initial set of variables,  $z\in\cit^n$, $t\in\cit^m$.

\begin{lemma}\label{l-rat2D'}
 For every $I\in P(\bm,n)$ we have
$$\prod_{(j,k,l,q)\in S_h(I)}\hspace{-15pt}(\phi^{(j,k,l,q)}_h-\phi^{(n+1,k,l,q)}_h)
\hspace{-15pt}\prod_{(j,k,l,q)\not\in S_h(I)}\hspace{-15pt}\phi^{(j,k,l,q)}_h=\sum_{I_0\in \mathcal{A}(h)} 
D'(F_{I_0}(h)(F_I(h)v)^*\otimes (F_{I_0}(h)v_{n+1})^*).$$
\end{lemma}
{\it Proof.} We will compare additive presentations of the two sides of the equality.   
The left hand side, after expansion of the products, has the following tree interpretation.
It is the rational function corresponding to a formal sum of forests. An element of this 
formal sum is a new forest produced from $T(I,h)$ as follows. We add
one tree, consisting only of strings of type $(l,k)$, $k\leq h<k$, corresponding to $\succ_h$
and grounded at zero. For each  string of type $(l,k)$ belonging to the new tree, one string of type
$(l,k)$ is removed from one of the trees which belong to the original forest $T(I,h)$. The coefficient
assigned to the new forest is a combinatorial coefficient times a sign coefficient, namely
minus one to the power equal to the number of strings which comprise the additional tree.  
Recall the definition of the index set 
$\mathcal{A}(h)=\{I_0=\{i_{l,k}\}_{k<l}\,|\,i_{l,k}=0 \mbox{ if } k>h, \mbox{ or } l\leq h\}$.
Thus, the  additional tree is always of the form $T(I_0,h,0)$ for some $I_0\in\mathcal{A}(h)$. 
In this notation the sign coefficient in the above description is $(-1)^{|I_0|}$. 
Let us give an explicit formula for the combinatorial coefficient. 
Consider $I_0\in\mathcal{A}(h)$, $I=(I_1,\ldots,I_n)\in P(\bm,n)$, 
where $I_j=\{i_{l,k}^{(j)}\}_{k<l}$ is a set of non-negative integers. In order to remove $i_{l,k}^0$ strings
of type $(l,k)$ from the forest $T(I,h)$, we should choose numbers $i_{l,k}^{(j,0)}$, where 
$0\leq i_{l,k}^{(j,0)}\leq i_{l,k}^{(j)}$ and $\sum_{j=1}^n i_{l,k}^{(j,0)}=i_{l,k}^{(0)}$.  Then, for each 
$j=1,\ldots,n$, we remove $i_{l,k}^{(j,0)}$ strings of type $(l,k)$ from the tree $T(I,h,z_j)\in T(I,h)$.
The combinatorial coefficient of this choice is 
$\prod_{k<l}\prod_{j=1}^n\bin{i_{k,l}^{(j)}}{i_{k,l}^{(j,0)}}$. Therefore, the left hand side equals
\begin{equation}\label{f-lhs}
\sum_{I_0\in\mathcal{A}(h)}(-1)^{|I_0|}\sum_{\{i_{l,k}^{(j,0)}\}}\left(\prod_{k<l}\prod_{j=1}^n
\bin{i_{k,l}^{(j)}}{i_{k,l}^{(j,0)}}\right)\phi(\{I-\sum_{j=1}^n\sum_{k<l}i_{k,l}^{(j,0)}1_{l,k}^{(j)},I_0\}, h).
\end{equation}

Next, we consider the right hand side. Repeated application of Lemma~\ref{F-act} gives
\begin{equation*}
 (-1)^{|I_0|} F_{I_0}(h)(F_I(h)v)^*=\sum_{\{i_{l,k}^{(j,0)}\}}\left(\prod_{k<l}\prod_{j=1}^n
\bin{i_{k,l}^{(j)}}{i_{k,l}^{(j,0)}}\right)
(F_{I-\sum_{j=1}^n\sum_{k<l}i_{k,l}^{(j,0)}1_{l,k}^{(j)}}(h)v)^*,
\end{equation*}
Therefore 
\begin{align}
&F_{I_0}(h)(F_I(h)v)^* \otimes (F_{I_0}(h)v_{n+1})^*= \label{eq-rhs}\\
& (-1)^{|I_0|}\sum_{\{i_{l,k}^{(j,0)}\}}\left(\prod_{k<l}\prod_{j=1}^n
\bin{i_{k,l}^{(j)}}{i_{k,l}^{(j,0)}}\right)(F_{I-\sum_{j=1}^n\sum_{k<l}i_{k,l}^{(j,0)}1_{l,k}^{(j)}}(h)v)^*
\otimes (F_{I_0}(h)v_{n+1})^*.\notag
\end{align}
Finally, formula~(\ref{f-lhs}), equality~(\ref{eq-rhs}), 
and the definition of the map $D'$ imply the statement of the lemma. $\qed$

\begin{lemma} \label{l-D'2D}
For every $I\in P(\bm,n)$, and every $I_0\in \mathcal{A}(h)$, we have
\begin{align}(a)\quad & D'(F_{I_0}(h)(F_I(h)v)^*\otimes 
(F_{I_0}(h)v_{n+1})^*)\backsim D(A\circ\tau(P_{I_0}(F(h),\Lambda_{n+1})) F_{I_0}(h)(F_I(h)v)^*).\notag\\
(b)\quad & D(A\circ\tau(P_{I_0}(F(h),\Lambda_{n+1})) F_{I_0}(h)(F_I(h)v)^*)\notag\\
 &=\sum_{J\in P(\bm,n)}\la (F_I(h)v)^*, A(F_{I_0}(h))\tau(P_{I_0}(F(h),\Lambda_{n+1}))F_J(h)v\ra D((F_J(h)v)^*).\notag
\end{align}
\end{lemma}
{\it Proof.} Part (a).
\begin{align}
 D' (F_{I_0}(h)(F_I(h)v)^*\otimes& 
(F_{I_0}(h)v_{n+1})^*)=  D'(F_{I_0}(h)(F_I(h)v)^*\otimes 
\tau(P_{I_0}(F(h),\Lambda_{n+1}))(v_{n+1})^*)\notag \\
& \backsim D'(A\circ\tau(P_{I_0}(F(h),\Lambda_{n+1}))
F_{I_0}(h)(F_I(h)v)^*\otimes (v_{n+1})^*)\notag \\
& =  D(A\circ\tau(P_{I_0}(F(h),\Lambda_{n+1}))
F_{I_0}(h)(F_I(h)v)^*).\notag
\end{align}
We use the property $(F_{I_0}(h)v_{n+1})^* = \tau(P_{I_0}(F(h),\Lambda_{n+1})(v_{n+1})^*)$ to write the first
equality, Corollary~\ref{switch-cor} to write the second equivalence, and the identification~(\ref{D'=D}) 
for the last step.

Part (b) is trivial. We write $A\circ\tau(P_{I_0}(F(h),\Lambda_{n+1})F_{I_0}(h)(F_I(h)v)^*$
in terms of the dual basis $\{(F_Jv)^*\}$ of $V[\nu]^*$, then use the linear property of the map $D$ and the 
definition of $sl_N$-action on the dual space.
\begin{align}
D&(A\circ\tau(P_{I_0}(F(h),\Lambda_{n+1})) F_{I_0}(h)(F_I(h)v)^*)=\notag\\
& = \sum_{J\in P(\bm,n)}D(\la A\circ\tau(P_{I_0}(F(h),\Lambda_{n+1}))
F_{I_0}(h)(F_I(h)v)^*,F_J(h)v\ra(F_J(h)v)^*)\notag\\
& = \sum_{J\in P(\bm,n)}\la A\circ\tau(P_{I_0}(F(h),\Lambda_{n+1})
F_{I_0}(h)(F_I(h)v)^*,F_J(h)v\ra D((F_J(h)v)^*)\notag\\
&= \sum_{J\in P(\bm,n)}\la (F_I(h)v)^*, A(F_{I_0}(h))\tau(P_{I_0}
(F(h),\Lambda_{n+1}))F_J(h)v\ra D((F_J(h)v)^*). \quad\qed\notag
\end{align}
 Finally, we combine all the steps in the following computation.
\begin{align}
&[\prod_{j=1}^n z_j^{(\Lambda_j,\omega_h^\vee)}
\prod_{d=1}^{m_h}(t_h^{(d)})^{-1}]\phi(z,t)= \notag\\
&=\prod_{j=1}^n z_j^{(\omega_h^\vee)^{(j)}}\hspace{-10pt}\sum_{I\in P(\bm,n)}
(\prod_{(j,k,l,q)\in S_h(I)}\hspace{-15pt}(\phi^{(j,k,l,q)}_h-\phi^{(n+1,k,l,q)}_h)
\hspace{-15pt}\prod_{(j,k,l,q)\not\in S_h(I)}\hspace{-15pt}\phi^{(j,k,l,q)}_h)F_Iv \notag\\
& = \prod_{j=1}^nz_j^{(\omega_h^\vee)^{(j)}}\hspace{-10pt}
\sum_{I\in P(\bm,n)}\sum_{I_0\in \mathcal{A}(h)}
D'(F_{I_0}(h)(F_I(h)v)^*\otimes (F_{I_0}(h)v_{n+1})^*)F_Iv\notag
\end{align}
\begin{align}
& = \prod_{j=1}^nz_j^{(\omega_h^\vee)^{(j)}}\hspace{-15pt}
\sum_{I,J\in P(\bm,n)}\sum_{I_0\in \mathcal{A}(h)}\la 
(F_I(h)v)^*,A(F_{I_0}(h))\tau(P_{I_0}(F(h),\Lambda_{n+1})(F_J(h))v\ra \times\notag\\
&\qquad\qquad\qquad \qquad\qquad\qquad 
\times D((F_J(h)v)^*)F_I(h)v\notag\\
& = \prod_{j=1}^nz_j^{(\omega_h^\vee)^{(j)}}\hspace{-10pt}
\sum_{J\in P(\bm,n)}\sum_{I_0\in \mathcal{A}(h)}D((F_J(h)v)^*)
A(F_{I_0}(h))\tau(P_{I_0}(F(h),\Lambda_{n+1})(F_J(h)v) \notag\\
&= \prod_{j=1}^nz_j^{(\omega_h^\vee)^{(j)}}\hspace{-10pt}
\sum_{J\in P(\bm,n)}D((F_J(h)v)^*)(\B_{\omega[h]}(\Lambda_{n+1}+\rho+\frac12 \nu)(F_J(h)v)).\notag
\end{align}
We have applied equality~(\ref{z-factor}), Lemma~\ref{l-rat2D'}, Lemma~\ref{l-D'2D}, 
contraction of the summation by $I$, and Theorem~\ref{main-add} consequtively.

Lemma~\ref{mat_trig_j1} asserts that $D((F_J(h)v)^*)=\phi(J,h)$. Since $\lambda=\Lambda_{n+1}+\rho+\frac12 \nu$
and $K_h(z,\lambda)=\prod_{j=1}^nz_j^{(\omega_h^\vee)^{(j)}}\B_{\omega[h]}(\lambda)$ we have
\begin{equation}
[\prod_{j=1}^n z_j^{(\Lambda_j,\omega_h^\vee)}
\prod_{d=1}^{m_h}(t_h^{(d)})^{-1}]\phi(z,t) \,\backsim\, K_h(z,\lambda)\phi(z,t),
\end{equation}
which is equivalent to the statement of the theorem for a hypergeometric solution $v(z,\lambda)$ 
with values in a tensor product of Verma modules. The statement of the theorem 
for a hypergeometric solution  with values in a tensor product of any highest weight $sl_N$-modules
follows from the functorial properties of the operator $K_h$.\hfill$\qed$

\section{APPENDIX A: Matsuo's hypergeometric solutions of the KZ equations for $sl_N$}
A construction  in \cite{Mat} gives hypergeometric solutions of the rational KZ equations in  
a tensor product of lowest weight $sl_N$ modules. We modify that procedure to 
a construction of solutions of the rational KZ equations in a tensor product of highest weight 
 $sl_N$ modules, cf. \cite{SV1}.

  Adopt the notation from Section~\ref{HGsol}. We have defined a  function
$\Phi':\cit^{n+1}_{z'}\times\cit^m_t\rightarrow\cit$, integration cycles $\gamma(z')$ 
and rational functions $\{\phi_I\}_{I\in P(\bm,n+1)}$. 

Order the basis $\{e_{N-1,N},\ldots,e_{1,2}\}$ of the Lie subalgebra $\gnp$ by
$e_{k,l}\succ e_{k',l'}$ if and only if $\alpha_{k,l}\succ\alpha_{k',l'}$. The corresponding
 PBW-basis of $U(\gnp)$ is
$$\left\{ E_{I_0} = \frac{e_{N-1,N}^{i_{N,N-1}}}{i_{N,N-1}!}\cdots
\frac{e_{1,2}^{i_{2,1}}}{i_{2,1}!}\right\}.$$
The index $I_0=(i_{l,k})_{l>k}$ runs over all sequences of non-negative integers.

 Let $W'=W_1\otimes\cdots \otimes W_{n+1}$, where $W_j$ is a lowest weight $sl_N$ module
of lowest weights $-\Lambda_j$ with lowest weight vector $w_j$. 
Fix a weight subspace
$W'[-\nu']\subset W'$, where $\nu'=\sum_{j=1}^{n+1}\Lambda_j-\sum_{i=1}^{N-1} m_i\alpha_i$ as in 
Section~\ref{HGsol}.

To every $I\in P(\bm,n+1)$, associate a vector 
$E_Iw=E_{I_1}w_1\otimes\cdots\otimes E_{I_{n+1}}w_{n+1}$. 

\begin{theorem}[Theorem~2.4 in \cite{Mat}]\label{t_mat_l} The function
$$w(z')=\int_{\gamma(z')}\Phi'(z',t)^{1/k}(\sum_{I\in P(\bm,n+1)}\phi_IE_Iw)\,dt$$
takes values in the subspace of singular vectors of $W'[-\nu']$
 and satisfies the rational KZ equations with parameter $\kappa\in\cit$.
\end{theorem} \hfill\\

\paragraph{\bf From lowest weight modules to highest weight modules.}
Recall that $\tau$ denotes the Chevalley involution. We have
$\tau(e_{k,l})=-e_{l,k}$ for any $1\leq l\ne k\leq N$. 

Let $\lambda\in\gh$, and let $V_{\lambda}$ be the  highest weight Verma module of 
highest weight $\lambda$ with highest weight vector $v_\lambda$, and let $W_{-\lambda}$ be the 
associated lowest weight Verma module of 
lowest weight $-\lambda$ with lowest weight vector $w_{-\lambda}$. The following proposition is well known.
\begin{proposition}\label{-to+}
The Chevalley involution defines an  isomorphism of Verma modules
$ W_{-\lambda}\rightarrow V_{\lambda}$ for generic $\lambda\in\gh$.
Namely, for all $x\in U(\gnp)$
we have $xw_{-\lambda}\mapsto\tau(x)v_\lambda$. This isomorphism induces  an isomorphism
from any lowest weight $sl_N$ module to the corresponding highest weight $sl_N$ module.
\end{proposition}  
 Assume that $(\Lambda_j)_{j=1}^{n+1}$ are generic.   
Let $V'=V_1\otimes\cdots\otimes V_{n+1}$, where $V_j$ is a highest weight 
module with highest weight $\Lambda_j$ and  highest weight vector $v_j$, 
corresponding to $W_j$ under the isomorphism of Proposition~\ref{-to+}. 
Notice that $\tau(E_J)=F_J$ for every set of positive integers $J=(j_{l,k})_{k<l}$. Thus,
$\sum_{I\in P(\bm,n+1)}\phi_IE_Iw\mapsto\phi'$. Theorem~\ref{t_mat_l} and Proposition~\ref{-to+} 
imply the following corollary.
\begin{corollary}\label{cor_mat}
The function
$$u(z')=\int_{\gamma(z')}\Phi'(z',t)^{1/k}\phi'(z',t)\,dt$$
takes values in the subspace of all singular vectors of $V'[\nu']$
 and satisfies the rational KZ equations with parameter $\kappa\in\cit$.
\end{corollary}

\section{APPENDIX B: Proof of Proposition~\ref{KZ2tKZ}}
Recall the statement of Proposition~\ref{KZ2tKZ}.\\
{\it
Fix a weight subspace $V'[\nu']\subset V'$, $\nu'=\sum_{j=1}^{n+1}\Lambda_j-\nu_0$, where  
$\bm\in Q_+$.  Let $u:\cit^{n+1}\rightarrow V'$ be a solution of the rational KZ 
equations with parameter $\kappa\in\cit$ taking values in the subspace of  $V'[\nu']$ 
consisting of all singular vectors. 
Set $\nu=\sum_{j=1}^{n}\Lambda_j-\nu_0$.

Then $v(z_1,\cdots,z_n)=u(z_1,\cdots,z_n,0)|v_{n+1}^*\ra
\prod_{i=1}^n z_i^\frac{(\Lambda_i,\Lambda_i+2\rho)}{2\kappa}$  is a solution of the 
trigonometric KZ equations with values in the weight subspace $V[\nu]\subset V$ with parameter
$\lambda = \Lambda_{n+1} +\rho +\frac12\nu\in\gh$ and the same parameter $\kappa\in\cit$.}

 {\it Proof.} 
 Set $u_0=u(z_1,\ldots,z_n,0)$. The function $u_0$ satisfies the system of equations
$\nabla_i(\kappa,\lambda)u_0=0$, $i=1,\ldots,n$.
Rewrite this equations  separating the input of the $(n+1)^{st}$ point in the sum:
\begin{equation}\label{eq1}
\kappa\frac{\partial u_0}{\partial z_i}= \left(\sum_{j=1, j\ne i}^n
 \frac{\Omega^{(ij)}}{z_i - z_j} +\frac{\Omega^{(i,n+1)}}{z_i}\right)u_0.
\end{equation}
Multiply equation~(\ref{eq1}) by $z_i$ and rearrange it using 
$\displaystyle\frac{z_i\Omega^{(ij)}}{z_i-z_j}=r(z_i/z_j)^{(ij)}+(\Omega^{-})^{(ij)}$
to obtain
\begin{equation*}
\kappa z_i\frac{\partial u_0}{\partial z_i}=\left(\sum_{j=1, j\ne i}^n
 \left(r(z_i/z_j)^{(ij)} + (\Omega^{-})^{(ij)}\right) + \Omega^{(i,n+1)}\right)u_0.
\end{equation*}
Set $v'=u_0|v_{n+1}^*\ra$. For all $\alpha\in\Sigma_+$, $i=1,\ldots,n$, we have
$e_{\alpha}^{(i)}v'=(e_{\alpha}^{(i)}u_0)|v_{n+1}^*\ra$.
$e_{-\alpha}^{(i)}v'=(e_{-\alpha}^{(i)}u_0)|v_{n+1}^*\ra$. This implies 
\begin{equation}
\kappa z_i\frac{\partial v'}{\partial z_i}=(\sum_{j=1, j\ne i}^n
 r(z_i/z_j)^{(ij)})v' + (\sum_{j=1, j\ne i}^n(\Omega^{-})^{(ij)}v' + 
(\Omega^{(i,n+1)})u_0|v_{n+1}^*\ra).\label{eqn6}
\end{equation}
 Now we claim that 
\begin{equation}\label{eqn7}
 (\Lambda_{n+1} +\rho +\frac12\nu-\frac12(\Lambda_i,\Lambda_i+2\rho))^{(i)}v'
=\sum_{j=1, j\ne i}^n(\Omega^{-})^{(ij)}v' + 
(\Omega^{(i,n+1)})u_0|v_{n+1}^*\ra).
\end{equation}
 First write 
\begin{equation}\label{eqn7.5}
(\Omega^{(i,n+1)})u_0|v_{n+1}^*\ra= (\sum_k x_k^{(i)}\otimes x_k^{(n+1)}u_0)|v_{n+1}^*\ra
+ \sum_{\alpha\in\Sigma_+}(e_{\alpha}^{(i)}\otimes e_{-\alpha}^{(n+1)}+
e_{-\alpha}^{(i)}\otimes e_{\alpha}^{(n+1)})u_0|v_{n+1}^*\ra.
\end{equation}

 The coupling $|v_{n+1}^*\ra$ allows us to compute the first summand explicitly
\begin{equation}\label{eqn8}
(\sum_k x_k^{(i)}\otimes x_k^{(n+1)}u_0)|v_{n+1}^*\ra = 
\sum_k(x_k,\Lambda_{n+1})x_k^{(i)}u_0)|v_{n+1}^*\ra=\Lambda_{n+1}^{(i)}u_0|v_{n+1}^*\ra=
\Lambda_{n+1}^{(i)}v'.
\end{equation}

Obviously $(e_{\alpha}^{(i)}\otimes e_{-\alpha}^{(n+1)})u_0|v_{n+1}^*\ra=0$.

Since $u_0$ is a singular vector we have 
$e_{\alpha}u_0=0$. Equivalently, $e_{\alpha}^{(n+1)}u_0=-\sum_{j=1}^ne_{\alpha}^{(j)}u_0$.
Therefore
\begin{align}\notag
&\sum_{\alpha\in\Sigma_+}(e_{-\alpha}^{(i)}\otimes e_{\alpha}^{(n+1)})u_0|v_{n+1}^*\ra=
\sum_{\alpha\in\Sigma_+}(e_{-\alpha}^{(i)}e_{\alpha}^{(n+1)})u_0|v_{n+1}^*\ra=
-\sum_{\alpha\in\Sigma_+}e_{-\alpha}^{(i)}(\sum_{j=1}^ne_{\alpha}^{(j)})u_0|v_{n+1}^*\ra\\
\notag &=-\sum_{\alpha\in\Sigma_+}(\sum_{j=1, j\ne i}^n 
e_{-\alpha}^{(i)}\otimes e_{\alpha}^{(j)})u_0|v_{n+1}^*\ra -
\sum_{\alpha\in\Sigma_+}(e_{-\alpha}^{(i)}e_{\alpha}^{(i)})u_0|v_{n+1}^*\ra
=-\sum_{j=1, j\ne i}^n(\Omega^{-})^{(ij)}u_0|v_{n+1}^*\ra -\\
\label{eqn9} &-(\frac12\sum_k x_k^{(i)} x_k^{(i)}+\sum_{\alpha\in\Sigma_+}
e_{-\alpha}^{(i)}e_{\alpha}^{(i)})u_0|v_{n+1}^*\ra+ 
\frac12\sum_{j=1}^n\sum_k x_k^{(i)} x_k^{(j)}u_0|v_{n+1}^*\ra.
\end{align}
Since $u_0\in V[\nu]$ we have 
$$\sum_{j=1}^n\sum_k x_k^{(i)} x_k^{(j)}u_0|v_{n+1}^*\ra=
\sum_k (x_k,\nu)x_k^{(i)}u_0|v_{n+1}^*\ra=\nu^{(i)}u_0|v_{n+1}^*\ra=\nu^{(i)}v'.$$
Let $C\in U(\gg)$ be the Casimir element. Since 
$e_{-\alpha}e_{\alpha}=-h_{\alpha}+e_{\alpha}e_{-\alpha}$ we have
$$-(\frac12\sum_k x_k^{(i)} x_k^{(i)}+\sum_{\alpha\in\Sigma_+}
e_{-\alpha}^{(i)}e_{\alpha}^{(i)})u_0|v_{n+1}^*\ra=-\frac12 C^{(i)}v' + 
\frac12\sum_{\alpha\in\Sigma_+}h_{\alpha}^{(i)} v'=-\frac12 C^{(i)}v' +\rho^{(i)} v'.$$
We have $C^{(i)}v'=(\Lambda_i,\Lambda_i+2\rho)v'$.
Rewrite (\ref{eqn9}) as
\begin{equation}\label{eqn10}
\sum_{\alpha\in\Sigma_+}(e_{-\alpha}^{(i)}\otimes e_{\alpha}^{(n+1)})u_0|v_{n+1}^*\ra=
(-\sum_{j=1, j\ne i}^n(\Omega^{-})^{(ij)}-\frac12 C^{(i)} +\rho^{(i)}+\frac12\nu^{(i)})v'.
\end{equation}
 Combine  (\ref{eqn10}),  (\ref{eqn8}) and
 (\ref{eqn7.5}) to obtain claim (\ref{eqn7}). Equation (\ref{eqn6}) for $v'$ becomes
\begin{equation}\label{eqn12}
\kappa z_i\frac{\partial v'}{\partial z_i}=\left(\sum_{j=1, j\ne i}^n
 r(z_i/z_j)^{(ij)}  +(\Lambda_{n+1} +\rho +
\frac12\nu)^{(i)}-\frac12\sum_{i=1}^n(\Lambda_i,\Lambda_i+2\rho)\right)v'.  
\end{equation}
Finally, $v=v'\prod_{i=1}^n z_i^\frac{(\Lambda_i,\Lambda_i+2\rho)}{2\kappa}$, and
equation~(\ref{eqn12}) implies  the proposition.\qed

\section{APPENDIX C: Formulae for $sl_3$.}
 Let $M_\lambda$ be the $sl_3$ Verma module with highest weight 
$\lambda\in\gh$ and highest weight vector $v_\lambda$.
Let $V$ be a tensor product of highest weight $sl_3$-modules.

Set $(\lambda,\alpha_j)=\lambda_j$, for $j=1,2$, \hspace{1cm}
$p_k(t)=t(t-1)\ldots(t-k+1)\in\cit[t]$ for any $k\in\nit$.

Formulae for the dynamical operators $\B_{\omega[1],V}(\lambda)=\B^{\alpha_1+\alpha_2}(\lambda)\B^{\alpha_1}(\lambda)$ 
and $\B_{\omega[2],V}(\lambda)=\B^{\alpha_1+\alpha_2}(\lambda)\B^{\alpha_2}(\lambda)$ are given in Section~\ref{add_sect}.

Let $\omega_0$ be the longest element of the Weyl group. It has two reduced presentations,
$\omega_0=s_1s_2s_1=s_2s_1s_2$, which imply two presentations  for the operator $\B_{\omega_0,V}$,\\
$\B_{\omega_0,V}(\lambda)=\B^{\alpha_1}(\lambda)\B^{\alpha_1+\alpha_2}(\lambda)\B^{\alpha_2}(\lambda)
=\B^{\alpha_2}(\lambda)\B^{\alpha_1+\alpha_2}(\lambda)\B^{\alpha_1}(\lambda)$. For any $v\in V[\nu]$, 
direct computation gives

\begin{align}
\B_{\omega_0,V}(\lambda+\rho+\frac12\nu)v&=
\sum_{a,b=0}^\infty\sum_{m,k=0}^{\min(a,b)}(-1)^mB_{m,k}^{a,b}(\lambda) e_{21}^{a-m}e_{31}^me_{32}^{b-m}
e_{12}^{a-k}e_{13}^{k}e_{23}^{b-k} v\notag\\
&=\sum_{a,b=0}^\infty\sum_{m,k=0}^{\min(a,b)}(-1)^mB_{m,k}^{a,b}(\lambda) e_{32}^{b-k}e_{31}^ke_{21}^{a-k}
e_{23}^{b-m}e_{13}^{m}e_{12}^{a-m} v,\notag
\end{align}
where, for any four non-negative integers $a,b,m,k$, such that $\max(m,k)\leq\min(a,b)$,
$B_{m,k}^{a,b}(\lambda)$ is a  complex-valued function depending on $\lambda\in\gh$
 defined by
$$B_{m,k}^{a,b}(\lambda)=\frac1{m!k!p_{a-m}(\lambda_1)p_{b-k}(\lambda_2)}\sum_{l=\max(m,k)}^{\min(a,b)}
\frac{(-1)^l l!}{(a-l)!(b-l)!(l-k)!(l-m)!p_l(\lambda_1+\lambda_1+1)}.$$

The shifted universal fusion matrix is
\begin{align}
 J_+(\lambda)&=\sum_{a,b=0}^\infty\sum_{m,k=0}^{\min(a,b)}
(-1)^{a+b}B_{m,k}^{a,b}(\lambda) e_{32}^{b-m}e_{31}^me_{21}^{a-m}\otimes e_{12}^{a-k}e_{13}^{k}e_{23}^{b-k}\notag\\
&=\sum_{a,b=0}^\infty\sum_{m,k=0}^{\min(a,b)}
(-1)^{a+b+m+k}B_{m,k}^{a,b}(\lambda) e_{21}^{a-k}e_{31}^ke_{32}^{b-k}\otimes e_{23}^{b-m}e_{13}^{m}e_{12}^{a-m}.\notag
\end{align}

The inverse to the Shapovalov form as an element of $M_\lambda\hat{\otimes} M_\lambda$ is
\begin{align}
 S_\lambda^{-1}&=\sum_{a,b=0}^\infty\sum_{m,k=0}^{\min(a,b)}
(-1)^{k}B_{m,k}^{a,b}(\lambda) e_{32}^{b-m}e_{31}^me_{21}^{a-m}v_\lambda\otimes 
e_{21}^{a-k}e_{31}^{k}e_{32}^{b-k}v_\lambda\notag\\
&=\sum_{a,b=0}^\infty\sum_{m,k=0}^{\min(a,b)}
(-1)^{k}B_{m,k}^{a,b}(\lambda) e_{21}^{a-k}e_{31}^ke_{32}^{b-k}v_\lambda\otimes 
e_{32}^{b-m}e_{31}^{m}e_{21}^{a-m}v_\lambda.\notag
\end{align}

Finaly, given $v\in V$, there is a unique singular vector, $sing(v_\lambda\otimes v)$, 
in $M_\lambda\otimes V$ of the form $sing(v_\lambda\otimes v)=  v_\lambda\otimes v+\,\{\mbox{lower order terms}\}$.
In \cite{ESt} the vector $sing(v_\lambda\otimes v)$ is given in terms of the inverse of 
the Shapovalov form. As a corollary we get
\begin{align}
 sing(v_\lambda\otimes v)&=\sum_{a,b=0}^\infty\sum_{m,k=0}^{\min(a,b)}
(-1)^{a+b}B_{m,k}^{a,b}(\lambda) e_{32}^{b-m}e_{31}^me_{21}^{a-m}v_\lambda\otimes 
e_{12}^{a-k}e_{13}^{k}e_{23}^{b-k}v\notag\\
&=\sum_{a,b=0}^\infty\sum_{m,k=0}^{\min(a,b)}
(-1)^{a+b+m+k}B_{m,k}^{a,b}(\lambda) e_{21}^{a-k}e_{31}^ke_{32}^{b-k}v_\lambda\otimes 
e_{23}^{b-m}e_{13}^{m}e_{12}^{a-m} v.\notag
\end{align}

\end{document}